\documentclass[11pt,a4paper,reqno]{amsart}
\usepackage{epsfig,amssymb,latexsym,amscd}

\usepackage[all,2cell]{xy}
\xyoption{v2}
\UseAllTwocells

\theoremstyle{plain}
\newtheorem{teo}{Theorem}[section]
\newtheorem{theo}[teo]{Theorem}
\newtheorem{coro}[teo]{Corollary}
\newtheorem{lema}[teo]{Lemma}

\theoremstyle{remark}

\theoremstyle{definition}

\newtheorem{defi}[teo]{Definition}
\newtheorem{obse}[teo]{Observation}

\newtheorem{example}[teo]{Example}

\newcommand{\id}{\ensuremath{\mathrm{id}}}

\newcommand{\su}{\operatorname{SU}_\mu(2,\mathbb{C})}

\begin{document}
\title[Compact quantum groups and the positive antipode]{Compact
  coalgebras, compact quantum groups and the positive antipode}
\author{Andr\'es Abella}
\thanks{The first author would like to thank Conycit-MEC, Uruguay.}
\address{Facultad de Ciencias\\Universidad de la Rep\'ublica\\
Igu\'a 4225\\11400 Montevideo\\Uruguay\\}
\email{andres@cmat.edu.uy}
\author{Walter Ferrer Santos}
\thanks{The second author would like to thank Csic-UDELAR,
  Conycit-MEC, Uruguay.}
\email{wrferrer@cmat.edu.uy}
\author{Mariana Haim}
\thanks{The third author would like to thank PEDECIBA, MEC-Udelar, Uruguay.}
\email{negra@cmat.edu.uy}

\begin{abstract} In this article --that has also the intention to survey
  some known results in the theory of compact quantum groups
using methods different from the standard and with a strong algebraic
  flavor-- we consider compact
$\circ$-coalgebras and Hopf algebras.
In the case of a $\circ$--Hopf algebra
we present a proof of
the characterization of the compactness in terms of the existence of a
positive definite integral, and use our methods to give an
elementary proof
of the uniqueness --up to conjugation by an automorphism of Hopf
algebras-- of the compact involution appearing in \cite{kn:andrus}.
We study the basic properties of
the positive square root of the antipode square that is
a Hopf algebra automorphism that we call the {\em positive
antipode}. We use it --as well as the {\em unitary antipode} and
Nakayama automorphism-- in order
to enhance our understanding of the antipode itself.
\end{abstract}
\maketitle
\section{Introduction}
\label{section:intro}

The family of {\em compact} groups plays an essential role
in the theory of group representations. This is clearly illustrated
in the pioneering work of  A. Hurwitz and I. Schur:
see for example Hurwitz's ``landmark
paper'' \cite{kn:hur} (c.f. Borel's \cite[Chapter 2, Section
2]{kn:borel}) or the later paper by I. Schur \cite{kn:schur}.

A specially illustrative example is the $unit\ddot{a}re\,\,
Beschr\ddot{a}nkung$ --later called unitarian trick-- that H. Weyl
used in order to
prove the complete reducibility of the representations of a semisimple
group --see \cite{kn:weyl1,kn:weyl2,kn:weyl3}--.

Since the mid 1980's it has been established
as a paradigm that the group objects in non commutative
geometry, \textit{i.e.} the so called {\em quantum groups}
are the --non necessarily commutative--
Hopf algebras.

Thus, after the above notion of {\em quantum group} was introduced,
it was natural  to expect the development of an adequate concept of
{\em compact quantum group}.

This concept was introduced in Woronowicz's seminal work --see
\cite{kn:woro}-- where the basic structure was defined in terms of a
$\star$--operation on a general Hopf algebra.

Later, the main {\em algebraic}
structure underlying the definition of compact quantum group
was introduced --and called  a {\em CQG algebra}--  where,
following the same spirit than in Woronowicz's paper, there is also
a compatible $\star$--operation  that
plays the central foundational role. See for example the paper
\cite{kn:andrus0} and the later \cite{kn:koor2} with a more formal
presentation. In book form, the reference \cite{kn:KS} has rather
complete a treatment of compact quantum groups and their representations.

In the present work we use the approach to the concept of compact
quantum group introduced in
\cite{kn:andres}, \cite{kn:andnic}, and \cite{kn:manin}, where
instead of the $\star$--operation
the equivalent concept of a $\circ$--structure is taken as the
fundamental notion. This approach is more
convenient for our purposes. In particular it allows the early
introduction of the basic concept of {\em compact coalgebra}.

In the following we present a brief description of the contents of
this paper.

In Section \ref{section:prelim}
we recall a few basic notations and well known results in coalgebra
and comodule theory that will be used throughout the paper.

In Section \ref{section:compactco} we recall the definition of
$\circ$--coalgebra and show that the existence of a
$\circ$--structure on a coalgebra $C$ is intimately
related with the existence of the analog of a {\em duality} functor in
the category of $C$--comodules. Then we introduce the notion of {\em compact
$\circ$--coalgebra} in terms of the existence of unitary inner products
in the objects of the category of $C$--comodules. In
particular, this definition shows
that compact $\circ$-coalgebras are cosemisimple.

In Section \ref{section:compactfourier} we give  different
characterizations of the notion of compact $\circ$-coalgebras in terms
of unitary inner products and also of positive definite Fourier
forms and present a
structure theorem that refines the corresponding well known result for
cosemisimple coalgebras --these results appear in
Theorem \ref{theo:coalgcompacta}--.

In Section \ref{section:hopf}
we specialize to the situation of Hopf algebras whose underlying coalgebra
and algebra structures have a compatible $\circ$--operator. In
Theorem \ref{theo:positivity}, we prove the following basic positivity
result: if $V$ and $V^*$ admit unitary inner
products, there exists a positive isomorphism between $V$ and
$V_{S^2}$ --where $V_{S^2}$ denotes de
$\mathbb{C}$--space $\,\/V$
with the structure corestricted with the automorphism $S^2$--.

In Section \ref{section:compactquantum} we define in terms of a
$\circ$--structure, the concept of compact quantum group and
present an intrinsic proof in Theorem \ref{theo:maincarac} --see also
Theorem \ref{theo:positivity2}-- of the
characterization of the compactness in terms of the positivity of the
bilinear  product defined from the integral --the Haar measure in
Woronowicz's nomenclature--. This proof is based upon the
existence --for every $H$--comodule $V$-- of a {\em positive}
isomorphism between $V$ and $V_{S^2}$ mentioned before.

Section \ref{section:compactinv} is dedicated to the presentation
of an intrinsic and
elementary algebraic proof of a
result due to Andruskiewitsch that appeared in \cite{kn:andrus}, and
that guarantees the
uniqueness up to Hopf algebra automorphisms, of the compact
involution, i,e., of the $\circ$--structure in the given Hopf
algebra.

In Section \ref{section:mien} we begin the study of some of
the basic properties of the
antipode in a compact quantum group. We show that for a
compact quantum group, there is a {\em positive} automorphism of $H$ that we
call the {\em positive antipode} --denoted as $S_+$-- that when
squared coincides with the automorphism $S^2$. Then, we use
some of the standard tools in coFrobenius Hopf algebra theory
--Radford's formula, Nakayama automorphism, the modular function,
etc.--
in order to
describe some of the basic properties of $S$ and of $S_+$. We show
that $S_+$ is given by conjugation by a positive multiplicative
functional $\beta$ --Theorem \ref{theo:S2inner} and Corollary
\ref{coro:squaremorf} --. We show that the Nakayama automorphism is a
{\em positive} algebra automorphism and define its postive square root
$\mathcal P$, that we write in terms of $\beta$. We compute explicitly
the adjoint of $S$ and show that the antipode is a normal operator
with respect to the inner product given by the integral,
if and only if the positive
antipode is the identity and this happens if and only if
$S^2=\operatorname {id}$ --Theorem
\ref{theo:meinanti}--. We also compute explicitly the so called {\em
  unitary antipode} in terms of $S$, $S_+$
and of the Nakayama automorphism.

We finish the paper with a short section --Section \ref{section:su}--
where  we illustrate the methods
developed, in the case of the compact quantum group
$\su$. In particular
we compute explicitly the positive antipode $S_+$ and the unitary
antipode $U$.

\medskip

The authors would like to thank Nicol\'as Andruskiewitsch who
--besides introducing them to this subject-- shared
many of his ideas, insights and widespread knowledge.

We also would like to thank the referee for her/his very useful
comments and corrections, that we are sure made the paper more
useful.

\section{Preliminaries and basic notations}
\label{section:prelim}
In this section we fix the notations and recall
a few  basic definitions and concepts in
coalgebra theory. See \cite{kn:andres, kn:koor} and
\cite{kn:take} for more details on some of the needed prerequisites
and \cite{kn:schneider,kn:sweedler}
for general background in Hopf algebra theory.

All the objects are defined over $\mathbb{C}$ and the
category of $\mathbb{C}$--vector
spaces will be denoted sometimes by $\mathcal V$.
If $V$ is a $\mathbb{C}$--vector space and
$A:V \rightarrow W$ is a linear map, then $V_c$ is the conjugate of $V$ and
$A_c: V_c \rightarrow W_c$ the conjugate linear map.
Moreover, $V^*_c$ will denote $(V^*)_c$.

For coalgebras and comodules we use Sweedler's notation.
The symbols: $\mathcal M^C$ and ${}^C\!\mathcal M$ represent the
categories of right and left $C$--comodules respectively and
$\mathcal M^C_f$ and ${}^C\!\mathcal M_f$
respectively are the full categories of finite dimensional objects.

\begin{defi} \label{defi:matrix} Let $V$ be a finite dimensional vector space and
  $\mathcal B= \{e_1,\cdots,e_n\}$ a basis.
Recall the definition of the coalgebra $c(V)=(V^* \otimes V, \Delta,
\varepsilon)$: $\Delta(\alpha \otimes v)= \sum \alpha \otimes e_i \otimes e^i \otimes
v$, and $\varepsilon: V^* \otimes
V \rightarrow \Bbbk$ is $\varepsilon(\alpha
\otimes v)= \alpha(v)$ --above $\mathcal B^*=\{e^1,\cdots,e^n\}$ is the dual
basis of $\mathcal B$--. The space $V$ is a right $c(V)$--comodule via
$\chi_0:V \rightarrow V \otimes c(V)$, given by $v \mapsto \sum e_i
  \otimes e^i \otimes v$. Then, $c(V)$ is a simple coalgebra and
$(V,\chi_0)$ is an irreducible right
$c(V)$--comodule (see for example  \cite[Section 3]{kn:andres}).
\end{defi}

\begin{obse} \label{obse:basic}
\begin{enumerate}
\item It is equivalent
to give a right $C$--comodule $(V,\chi) \in \mathcal M^C_f$ and
to give a morphism of coalgebras $c_{\chi}: c(V)
  \rightarrow C$. The morphism of coalgebras
is the map $c_{\chi} (\alpha \otimes v)= \sum
  \alpha(v_0)v_1$.
Notice that in the case that $V$ is a not necessarily finite
  dimensional $C$--comodule, we can still define the linear map
  $c_{\chi}: V^* \otimes V \rightarrow C$ but
the vector space $V^* \otimes V$ is not endowed with a
  natural coalgebra structure but $c_{\chi}(V^* \otimes V) = \operatorname{Coeff}(V)\subset C$ is
  a subcoalgebra of $C$, called the \emph{subcoalgebra of
    coefficients} of $V$.
In that situation the $C$--comodule structure of $V$ is induced from a
$\operatorname{Coeff}(V)$--comodule structure. In other words,
$\chi: V \rightarrow V \otimes C$ can be split as $V
  \xrightarrow{\chi} V \otimes \operatorname{Coeff}(V) \subset V
  \otimes C$.
In this situation the
relationship between the $C$--comodule structure of $V$ and the maps $\chi_0$
and $c_\chi$ introduced respectively in Definition \ref{defi:matrix} and Observation
\ref{obse:basic} part (1), is given by
the commutativity of the diagram below.

\[\xymatrix{V \ar[d]_{\chi_{_0}} \ar[rr]^-{\chi}& & V \otimes C
\\ V \otimes c(V) \ar[rr]_-{\operatorname {id} \otimes c_{\chi}}
 & & V \otimes \operatorname{Coeff}(V) \ar[u]_{\operatorname{id}
   \otimes \operatorname{inc}}}\]
\item In the above notations, given $(V,\chi) \in \mathcal M^C_f$,
we define the matrix $T \in \operatorname
  {M_n}(C)$, as $T=(c_{\chi}(e^j \otimes e_i))_{1\leq i,j \leq n}$.
The {\em matricial
coefficients} of the $C$--comodule $V$, i.e. the entries of $T$,
are defined by the formul\ae\/\, $\chi(e_j)=\sum_i e_i \otimes t_{ij}$.
Clearly, $\Delta(t_{ij})=\sum_k t_{ik}
  \otimes t_{kj}$\,; $\varepsilon(t_{ij})=\delta_{ij}$.  It is well
  known that $V$ is irreducible if and only if the matrix
  coefficients $\{t_{ij}\}$ are linearly independent.
\end{enumerate}
\end{obse}

In the context of linear actions of groups, the
 algebra of coefficients is usually called the {\em algebra of
 representative functions}.

We list a few results concerning the above construction. Most of the
proofs  can be found in \cite[Section 4]{kn:andres}.
\begin{obse} \label{obse:basiccoeff}
 \begin{enumerate}
\item Assume that $I \subset C$ is a right coideal with the property
  that $\operatorname{Coeff}(I) \subset I$. Then $\Delta(I) \subset I
  \otimes \operatorname{Coeff}(I)$ and applying $\varepsilon \otimes
  \operatorname {id}$ we deduce that $I \subset
  \operatorname{Coeff}(I)$. Hence $I=\operatorname{Coeff}(I)$ is a
  subcoalgebra.
\item If $V \in
\mathcal M^C$ is irreducible, then $\operatorname{Coeff}(V)$ is a
simple subcoalgebra of $C$ and conversely if $C$ is
a simple coalgebra, there exists an irreducible $C$--comodule $V$ such
that $C=\operatorname {Coeff}(V) \cong c(V)$. Moreover if $C$ is
a simple and $W$ is an arbitrary non zero
$C$--comodule, then $\operatorname {Coeff}(W)=C$.
\item If $C$
is an arbitrary coalgebra and $C_1= \operatorname {Coeff}(V_1),\,
C_2= \operatorname {Coeff}(V_2)\subset C$ are simple
subcoalgebras  with $V_1,V_2$ irreducible $C$-comodules, then $C_1=C_2$ if and only if
$V_1 \cong V_2$ as $C$--comodules.
\item If $C$ can be decomposed into a sum of irreducible
  $C$--comodules: $C=\bigoplus \{V_\rho: \rho \in \mathcal
  R\}$, then an arbitrary $C$--comodule $V$, is isomorphic to
  $V_\rho$ for some $\rho \in \mathcal R$ and
any simple subcoalgebra $D \subset C$ is of the form
$D=\operatorname {Coeff}(V_\rho)$ for some $\rho \in \mathcal R$.
\end{enumerate}
\end{obse}

\begin{obse}\label{obse:positivity}

A hermitian form $\beta: V \otimes V_c \rightarrow \mathbb{C}$
defined on $V$ induces a linear morphism
$\psi_\beta:V \rightarrow V^*_c$ defined as
$\psi_\beta(u)(v)=\beta(v,u)$.

Assume that  $V$ and $V^*$ are finite dimensional and
endowed with inner products $\beta: V \otimes V_c \rightarrow
\mathbb{C}$
and $\gamma: V^* \otimes V^{*}_c
\rightarrow \mathbb{C}$. Next, consider the corresponding
isomorphisms $\psi_\beta:V \rightarrow V^{*}_c$ and
$\psi_\gamma:V^*\rightarrow V^{**}_c$ and
the linear map $\phi:V \rightarrow V$ defined by the
commutativity of the following diagram:
\[
\xymatrix{
V \ar[d]_{{\psi_\beta}} \ar[r]^-{\phi} & V  \ar[d]^{j_V} \\
V^*_c \ar[r]_-{{\psi_{\gamma_c}}}          & V^{**}
}
\]
where $j_V: V \rightarrow V^{**}$ denotes the canonical isomorphism,
i.e. $\phi$
is defined by the equality: $\gamma(
\psi_\beta(v),\psi_\beta(u))=\beta(\phi(u), v),\,\, \forall u,v\in V$.

Then $\beta(\phi(v), v)>0$, for all $0\ne v\in V$, i.e.
$\phi:V \rightarrow V$ is $\beta$--positive definite.
\end{obse}

Some aspects of the theory of $\circ$--coalgebras, unitary
comodules and compactness --that will be treated later--
can be viewed with profit in terms of Fourier forms and Fourier
products in the coalgebra. The theory of Fourier
forms and Fourier products for general coalgebras was developped in
\cite{kn:ferrerfourier}, we recall here a few basic definitions, the
fact that the basic field is $\mathbb{C}$ is not relevant here.

\begin{defi} \label{defi:transform}
\begin{enumerate}
\item A \emph{Fourier form} in $C$ is a bilinear map $\omega: C \otimes C
  \rightarrow \mathbb{C}$ with the property that  for all $c,d \in C$,
  $\sum \omega(c, d_1)d_2=\sum c_1 \omega(c_2, d)$. The
  Fourier form is said to be \emph{normal} if for all $c \in C$,
  $\sum\omega(c_1, c_2)=\varepsilon(c)$. The Fourier form is
  \emph{symmetric} if $\omega(c, d)= \omega(d, c)$ for all $c,d
  \in C$.
\item A \emph{Fourier product} in $C$ is a bilinear map $\star: C \otimes C
  \rightarrow C$ with the property that  for all $c,d \in C$,
  $\Delta(c\star d)=\sum (c \star d_1)\otimes d_2=\sum c_1 \otimes
  (c_2 \star d)$.

\item If $H$ is a Hopf algebra and $\varphi:H \rightarrow \mathbb{C}$ is a
  right integral on $H$ \textit{i.e.} a linear map such that 
$\sum\varphi(x_1)x_2=\varphi(x)1$, for all $x\in H$\footnote{In a first version of this paper and
    also in some previous papers by the authors --see for example
    \cite{kn:ferrerfourier}-- the name cointegral was used for this
    kind of maps. We changed it following the referee's suggestion. The present
    name is more compatible with the usual nomenclature on the subject.},
then $\omega_{\varphi}(x, y)=\varphi(S(x)y)$ is a Fourier form --see \cite{kn:ferrerfourier}--.
Moreover, the integral $\varphi$ is normal --i.e. it satisfies that
$\varphi(1)=1$-- if
and only if $\omega_{\varphi}$ is a normal Fourier form. Conversely, if
$\omega$ is a Fourier form, then the linear map
$x \mapsto \omega(1, x): H \rightarrow \mathbb{C}$ is a right integral
on $H$.

\end{enumerate}
\end{defi}

Below we establish a bijective correspondence between Fourier forms
and Fourier
products on a coalgebra $C$.
Indeed:
$\omega \mapsto \star_{\omega}$, $\star \mapsto \omega_\star$,
where for all $c,d \in C$,
\[
c\star_\omega d= \sum \omega(c , d_1)d_2=\sum c_1 \omega(c_2 , d)\,\,
, \,\, \omega_\star( c,d) = \varepsilon(c \star d)\,\, ,
\]

It is clear that the Fourier form $\omega$ is
normal if and only if the associated Fourier product satisfies
the following condition: $\sum x_1 \star_\omega x_2= x$.

The following properties are easy to verify.

\begin{obse}\label{obse:lem1}
\begin{enumerate}

\item If $\star_1$ and $\star_2$ are arbitrary Fourier products in $C$,
then $x \star_1 (y\star_2 z)= (x\star_1 y) \star_2 z$.

\item If $\star$ is a Fourier product in $C$ and $L,R$ are respectively
a left and a right coideal in $C$, then $R\star L\subseteq R\cap
L$. Hence, if $D$ and $E$ are two different simple subcoalgebras of
$C$, then $D\star E=0$ and therefore $\omega_\star(D\otimes E)=0$.

\item If $C$ is cosemisimple and
$\omega:C\otimes C \rightarrow \mathbb{C}$ is a right non degenerate
Fourier form, then
it follows from the above that
for every simple subcoalgebra $D \subseteq C$, the restricted
form $\omega:D \otimes D\rightarrow \mathbb{C}$ is right (and left --since
$D$ is finite dimensional--) non degenerate.

\item If $D$ is a finite dimensional coalgebra and $\omega, \star$
are a form and a product associated to each other as above, then: $\omega:D\otimes D \rightarrow \mathbb C$ is right
(and left) non degenerate if and only if
$\star$ has a neutral element $s\in D$.
Explicitely,  $s$ is the only element in
$D$ that verifies
$\omega (c, s)=\varepsilon(c)=\omega(s, c), \forall c \in D$.

\end{enumerate}
\end{obse}

\section{Compact $\circ$--coalgebras}
\label{section:compactco}

If $G$ is a compact topological group and $\mu$ is an invariant
measure, a standard tool in the representation theory of $G$
consists in ``unitarizing'' all the $G$--modules.
Given an arbitrary inner product in a $G$--module
$V$, by integration we define on $V$ a $G$--invariant inner
product and in that manner we can assume that the action of the group on
$V$ is given by {\em isometries}.

In order to generalize
the theory of compact groups to general coalgebras, it is natural to
assume as the starting point of the abstract definition of compactness
of $C$, the existence of an  inner product on every comodule
satisfying a condition of $C$--invariance --see more precisely Definition
\ref{defi:unibil}--.

This will lead immediately to the concept of
compact $\circ$--coalgebra that is presented in this section.

The notion of a $\circ$--structure in the context of {\em Hopf algebras}
was originally defined in \cite{kn:manin}, where the concept was
presented in terms of a
$\star$--structure and of the antipode.

Later the definition of a
$\circ$--structure in the more general set up of {\em coalgebra theory}
appeared in  \cite{kn:andnic} (see also
\cite{kn:andres} for a more recent and detailed
presentation).

The algebra dual to a $\circ$--coalgebra is a $\star$--algebra, where
the operator $\star$ is the dualization of $\circ$. Many concepts
and proofs in the realm of $\circ$--coalgebras are counterparts of
concepts and proofs for $\star$--algebras.

\begin{defi}\label{defi:cerocoalgebra}
\begin{enumerate} 
\item 
Assume that $(C,\Delta,\varepsilon)$ is a coalgebra over the field of complex numbers. 
A $\circ$--{\em structure} on $C$ is a linear map $\circ: C\rightarrow C_c$ that is involutive and 
anticomultiplicative. The pair $(C,\circ)$ is called a \emph{$\circ$--coalgebra}.
\item 
If $(C,\circ)$ and $(D,{\lozenge})$ are two $\circ$--coalgebras, a
\emph{morphism} of $\circ$--coalgebras is a morphism of coalgebras $f:C\rightarrow D$, 
satisfying the additional condition  that for all $c\in C$, $f(c^{\circ})=f(c)^{{\lozenge}}$.
\item 
A subcoalgebra $D$ of a $\circ$--coalgebra satisfying
$D^\circ = D$ is called a $\circ$--\emph{subcoalgebra}.
\end{enumerate}
\end{defi}

In explicit terms and with respect to the comultiplication,
the compatibility condition for $\circ$ means that for all $c \in C$ if
$\Delta(c)=\sum c_1 \otimes c_2$, then $\Delta(c^{\circ})= \sum
c^{\circ}_2 \otimes c^{\circ}_1$.

\begin{obse}\label{obse:epsilon}
It follows immediately from the
  definition above that
  $\varepsilon(c^{\circ})=\overline{\varepsilon(c)}$ for all $c \in
  C$. Hence viewing $(C_c,\Delta_c,\overline{\varepsilon_c})$ as a
coalgebra, the $\circ$--structure is a morphism of coalgebras
$\circ:(C,\Delta,\varepsilon) \rightarrow
(C_c,\Delta_c,\overline{\varepsilon_c})^{\mathrm{cop}}$.
\end{obse}

\begin{example} \label{example:basic}

If $V$ is an arbitrary finite dimensional complex vector space with
basis $\mathcal B=\{e_1,\cdots,e_n\}$, then the coalgebra $c(V)= V^* \otimes
V$ can be equipped with a $\circ$--structure 
defined on the basis $\mathcal B$ as $(e^j \otimes e_i)^\circ=e^i \otimes e_j$.
\end{example}

The $\circ$--structure on the coalgebra $C$ endows the category of
right $C$--comodules with an additional ``duality functor'' as presented
in Definition \ref{defi:dual}. The reader should notice that we
are not talking about duality in a technical sense.
\medskip

First we recall the following duality constructions that make sense in the
context of general coalgebras.

\begin{defi}\label{defi:radjoint}
 Let $C$ be an arbitrary coalgebra and $(V,\chi)$ a right
 finite dimensional $C$--comodule. The map $\chi^r: V^*
  \rightarrow C \otimes V^*$ defined by the formula:
$\chi^r(f)=\sum f_{-1} \otimes f_0$ if and only
  if $\sum f_{-1}f_0(v)= \sum f(v_0)v_1$, is a left $C$--comodule
  structure on $V^*$. The pair $(V^*,\chi^r)$ is
  frequently abbreviated as $V^r$ and called the \emph{right adjoint} of
  $V$. In this manner we have defined a contravariant functor
$(-)^r:\mathcal M^C_f
  \rightarrow {}^C\!\mathcal M_f$. Similarly we can define
  $(-)^{\ell}: {}^C\!\mathcal M_f \rightarrow \mathcal M^C_f$ and
  $V^{\ell}$ will
  be called the \emph{left adjoint} of $V$.
\end{defi}

\begin{defi}
If $C$ is a $\circ$--coalgebra and $(W,\chi)$
a left $C$--comodule, then
the map $\chi^\circ$
defined as $\chi^\circ=\operatorname{sw}(\circ \otimes \operatorname{id})\chi_c:
W_c \rightarrow W_c \otimes C$ is a right $C$--comodule structure on
$W_c$, where $\operatorname{sw}$ stands for the usual switching map.
This correspondence can be extended to define
a covariant functor ${}^C\!\mathcal M\to \mathcal M^C$.
\end{defi}

\begin{defi}\label{defi:dual}
Let $C$ be a $\circ$--coalgebra. If $(V,\chi) \in \mathcal M^C_f$, then the map
$\chi^{r\circ}=\left(\chi^r\right)^\circ :V^*_c\rightarrow V^*_c \otimes C$
is a right $C$--comodule structure on $V^*_c$.
In explicit terms, if $\chi(v)=\sum v_0\otimes v_1$ and
$\chi^{r\circ}(f)=\sum f_0\otimes f_1$ then
\[
\sum \overline{f_0(v)} f_1=\sum \overline{f\left(v_0\right)} v_1^\circ.
\]
We call $\mathcal D$ the contravariant functor
$\mathcal M^C_f \rightarrow \mathcal M^C_f$ defined on objects as
above, i.e.
$\mathcal D(V,\chi)= (V^*_c,\chi^{r\circ})$.
If $\mathcal B=\{e_1,\cdots,e_n\}$ is a basis of $V$ and $\mathcal B^*=\{e^1,\cdots,e^n\}$ is
its dual basis in $V^*$, then the matrix coefficients of $\chi$ and $\chi^{r\circ}$ are related by
\[
\chi(e_i)=\sum_j e_j\otimes t_{ji},\quad \chi^{r\circ}\!\left( e^i \right)=\sum_j e^j\otimes t_{ij}^\circ
\]
\end{defi}

\begin{obse}\label{obse:coeffdual}
In the situation above, it is clear that for a finite dimensional
$C$--comodule $V$, we have that
$\operatorname{Coeff}(\mathcal D(V))= \operatorname{Coeff}(V)^\circ$.
\end{obse}

Next lemma follows immediately from the fact that $\circ$ is an involution.

\begin{lema}\label{lema:ddual}
If $C$ is a $\circ$--coalgebra, then
$\mathcal D^2=(-)^{{\ell}r}: \mathcal M^C_f \rightarrow \mathcal M^C_f$.
\qed
\end{lema}

In order to look at the adequate concept of representation for
a $\circ$--coalgebra we need the following definition.

\begin{defi}\label{defi:unibil}\label{defi:unitary}
\begin{enumerate}
\item Let $C$ be a $\circ$--coalgebra and $V \in\mathcal M^C$.
A hermitian form $\beta : V \otimes V_c \rightarrow \mathbb{C}$ is
\emph{unitary} or \emph{invariant} if $\sum \beta(u_0,v)u_1 = \sum
\beta(u,v_0) v_1^{\circ}$, for all $u,v \in V$. If moreover,
$\beta$ is an inner product, then it is called a \emph{unitary
  structure} on $V$.
\item A pair $(V,\langle\,\,\,,\,\,\rangle)$ consisting of a $C$--comodule $V$ and a unitary
  structure on $V$ is called a \emph{unitary} $C$--comodule.
Moreover, if $(V,\langle\,\,\,,\,\,\rangle)$ and $(W,[\,\,\,,\,\,])$
  are unitary comodules, a morphism of comodules $f:V \rightarrow W$ is
  said to be \emph{unitary} if it is an isometry with respect to
  $\langle\,\,\,,\,\,\rangle$ and $[\,\,\,,\,\,]$. Sometimes a
  \emph{unitary structure} is called a \emph{unitary inner product}.
\end{enumerate}
\end{defi}

\begin{example}
Let $V$ be a finite dimensional vector space and consider the
$\circ$--coalgebra $c(V)$ of Example \ref{example:basic}. Then, the
inner product of
$V$ that makes the given basis $\mathcal B$
orthonormal, is a unitary form in the $c(V)$--comodule $V$.
\end{example}

\begin{lema}\label{lema:unitary}
Let $C$ be a $\circ$--coalgebra, $(V,\chi)\in\mathcal M^C_f$ and let
$\beta: V \otimes V_c \rightarrow \mathbb{C}$ be a hermitian form on
$V$.

\noindent
Then, the following assertions are equivalent:

\medskip
\noindent (1) The form $\beta$ is unitary.\\
(2) The map $\psi_\beta:V = (V,\chi) \rightarrow \mathcal D(V)= (V^*_c,\chi^{r \circ})$
--see Observation \ref{obse:positivity}-- is a morphism in $\mathcal
M^C_f$.
\end{lema}

\begin{coro}\label{coro:nondeg}
If $\beta$ is a unitary structure on the $C$--comodule $V$, then
$\operatorname{Coeff}(V) \subset C$ is a $\circ$--subcoalgebra that is
simple if $V$ is irreducible.
\end{coro}
\begin{proof}
In this situation $\psi_\beta$ is an isomorphism,  then
the result follows from Observations \ref{obse:basiccoeff}
and \ref{obse:coeffdual} from which we deduce that
$\operatorname{Coeff}(V)=\operatorname{Coeff}(V^*_c)=\operatorname{Coeff}(V)^{\circ}$.
\end{proof}

\begin{coro}\label{coro:nonisoorth}

Assume that $C$ is a $\circ$--coalgebra, that $V$ is a $C$--comodule and that
$\beta: V \otimes V_c \rightarrow \mathbb{C}$ is a
  unitary structure on $V$. If $V_\lambda$ and $V_\mu$ are
  irreducible non isomorphic subcomodules of $V$, then $V_\lambda$ and
$V_\mu$ are $\beta$--orthogonal.
\end{coro}
\begin{proof} Consider $C_\lambda=\mathrm{Coeff}(V_\lambda) \subset C
  \,\,\/\text{and}\,\,C_\mu=\mathrm{Coeff}(V_\mu)
  \subset C$, the simple subcoalgebras of coefficients associated
  to $V_\lambda$ and
  $V_\mu$ respectively. As $V_\lambda$ and $V_\mu$ are not isomorphic,
  then $C_\lambda \cap C_\mu=\{0\}$ --see Observation
  \ref{obse:basiccoeff}--.
In this situation, for $u \in V_\lambda$ and $v \in V_\mu$, if we apply
  $\varepsilon$ to the equality  $\sum \beta(u_0, v) u_1 = \sum
  \beta(u, v_0) v_1^{\circ} \in C_\lambda \cap
  C_\mu ^{\circ}= C_\lambda \cap
  C_\mu =\{0\}$, we deduce that $\beta(u, v)=0$.
\end{proof}

In the presence of a unitary structure the comments appearing in
Observation \ref{obse:basic} can be refined.

\begin{lema} \label{lema:unit_trans}
Assume that $C$ is an arbitrary
  $\circ$--coalgebra and that $(V,\chi)$ is a finite dimensional right
  $C$--comodule. Choose a basis $\mathcal B=\{e_1,\cdots,e_n\}$ and
  define the matrix  $T=(t_{ij})_{1 \leq
  i,j \leq n}$ by the formula $\chi(e_i)=\sum_j e_j \otimes t_{ji}$.
 Then $T^\circ=T^t$ if and only if the inner product that makes
 the basis $\mathcal B$ orthonormal, is a unitary structure.
\end{lema}
\begin{proof}
First we observe that if the product is unitary, then
the condition $T^\circ=T^t$ is satisfied.
The unitary condition $\sum \langle u_0,v\rangle u_1 =
\sum \langle u, v_0\rangle v_1^\circ$ applied to a pair of basis
vectors $u=e_k$ and $v = e_j$ yields directly the
equality $t_{ki}^\circ = t_{ik}$, \textit{i.e.},
$T^\circ=T^t$. Conversely, if this condition is satisfied, the above
reasoning can be reversed and one proves that the inner product is unitary.
\end{proof}

\begin{obse}\label{obse:unit_cr}
If $C$ is a $\circ$--coalgebra and
$(V,\langle\,\,\,,\,\,\rangle)$ is
  a unitary $C$--comodule, it is clear from the definitions
that if $W \subset V$
  is an arbitrary $C$--subcomodule,
then the orthogonal complement
$W^{\perp}$ of $W$ is also a
  $C$--subcomodule of $V$ (see also \cite[Corollary 5.16]{kn:andres}
  for details).
Hence a unitary $C$--comodule is completely
  reducible.
\end{obse}

\begin{defi} \label{defi:compactcoalgebra}
If $C$ is a $\circ$--coalgebra, we say that $C$ is
\emph{compact} if every right $C$--comodule admits a unitary structure as
introduced in Definition \ref{defi:unitary}.
\end{defi}

\begin{obse} \label{obse:muchascosas} \begin{enumerate}
\item It follows from Observation \ref{obse:unit_cr} that a
  compact coalgebra is co\-semi\-simple.
\item If $C$ is compact, and we consider it as a right $C$--comodule,
  then it can
  be endowed with a unitary inner
product that will satisfy the following condition: for all $c,d \in
C$, $\sum \langle c_1,d \rangle c_2 = \sum \langle c,d_1 \rangle d_2^{\circ}$.

\item In the situation above if $1 \in \mathbb{C}$ is a group like
  element with the property that $1^{\circ}= 1$, we can define a map $\varphi: C
  \rightarrow \mathbb{C}$, as $\varphi(c)=\langle c,1
  \rangle \in \mathbb{C}$. It is clear that $\varphi$ is a right
  integral on $C$ with $\varphi(1)=1$. Indeed, $\sum \varphi(c_1)c_2= \sum \langle c_1,1
  \rangle c_2= \langle c,1 \rangle 1^{\circ}= \varphi(c)1$. Moreover
  $\varphi(1)= \langle 1,1 \rangle >0$ because we are dealing with inner
  products.
\end{enumerate}
\end{obse}

\section{Compactness, cosemisimplicity and Fourier forms}
\label{section:compactfourier}

In the presence of a $\circ$--structure on $C$ the existence of inner
products can be related with the theory of Fourier forms and some of
the considerations of the preceeding section can be refined.

\begin{defi}\label{defi:positive} Assume that $C$ is a
  $\circ$--coalgebra, that $\omega$ is a Fourier form for $C$ and that
  $\star$ is a Fourier product.
\begin{enumerate}
\item We say that $\omega$ is {\em positive} if and only if
  $\omega(c^\circ, c) \geq 0$ for all $c \in C$. If $\omega(c^\circ, c) > 0$ whenever $c \neq 0$ we say that $\omega$ is
  {\em positive definite}.
\item  We say that $\omega$ is {\em hermitian} if for all $c,d \in C$,
  $\omega\!\left(c^\circ, d^\circ\right)=
  \overline{\omega(d, c)}$.
\item We say that $\star$ is {\em positive} if and only if
  $\varepsilon(c^\circ \star c) \geq 0$ for all $c
  \in C$. If $\varepsilon(c^\circ
  \star c) > 0$ whenever $c \neq 0$ we say that $\star$ is
  {\em positive definite}.
\item We say that $\star$ is {\em hermitian} if for all $c,d \in C$,
  $c^\circ \star\, d^\circ = (d\star c)^\circ$.
\end{enumerate}
\end{defi}

It follows directly from the above definition that the Fourier form
$\omega$ is positive (positive definite, hermitian) if and only if the
corresponding Fourier product $\star_\omega$ is positive
(respectively positive definite, hermitian).

\begin{obse}\label{obse:refinamiento}
\begin{enumerate}
\item In the case that $C$ has a $\circ$--structure and also a Fourier form
$\omega$, the map  $\langle\,\,\,,\,\,\rangle_\omega: C \otimes C \rightarrow
\mathbb{C}$ defined for all $c,d \in C$ as $\langle c,d
\rangle_\omega= \omega(d^\circ, c)$ is a unitary sesquilinear form in $C$.
\item Hence, in the situation of a $\circ$--coalgebra $C$,
there is a bijective correspondence between
  Fourier forms, Fourier products  and unitary sesquilinear forms on
  $C$ and this
  correspondence is given by the following rules:  $\omega \mapsto
  \star_{\omega}$, $\star \mapsto \langle\,\,\,,\,\, \rangle_\star$
  and $\langle\,\,\,,\,\,\rangle \mapsto
  \omega_{\langle\,\,\,,\,\,\,\rangle}$ where for all $c,d \in C$,
$ c\star_\omega d=
  \sum \omega(c, d_1)d_2 = \sum c_1 \omega(c_2, d)$, $\langle c, d\rangle_\star = \varepsilon(d^\circ \star c)$ and
$\omega_{\langle\,\,\,,\,\,\,\rangle}(c, d) = \langle d,
c^\circ \rangle$.
\item The direct expression of the Fourier product in terms of the
  corresponding unitary sesquilinear form is the following: for all
$c,d \in C$: $c \star d = \sum \langle d_1,c^\circ\rangle d_2=\sum c_1\langle
  d,c_2^\circ\rangle$.
\item The Fourier form $\omega$ is symmetric if and only if the
  associated sesquilinear form satisfies the following condition:
  $\langle c^\circ,d^\circ\rangle= \langle d,c \rangle$.
In case that the sesquilinear form is hermitian, the symmetry of
  $\omega$ is equivalent to the condition $\langle
  c^\circ,d^\circ\rangle=\overline{\langle c,d \rangle}$.
\item  The correspondence considered above between
  Fourier forms and unitary sesquilinear forms in $C$ , preserves positivity
  and hermitianity.
  Hence, in $C$ there is
  a bijective correspondence between positive definite hermitian
  Fourier forms and unitary inner products.
\end{enumerate}
\end{obse}

The theorem below refines for the case of compact coalgebras, a well
known result on the structure of cosemisimple coalgebras.

\begin{theo}\label{theo:coalgcompacta}
Let $C$ be a $\circ$--coalgebra. The following properties are equivalent.
\begin{enumerate}
\item \label{e:compactcoalg} The coalgebra $C$ is compact.
\item \label{e:unitinnerpro} The coalgebra $C$ viewed as a right $C$--comodule admits a
  unitary inner product.
\item \label{e:fourierpos}
The coalgebra $C$ admits a Fourier form that is positive definite and
  hermitian.
\item \label{e:prodpos}
The coalgebra $C$ admits a Fourier product that is positive definite and
  hermitian.
\item \label{e:decomposition}
The coalgebra $C$ can be decomposed as
$C=\bigoplus_{\rho \in \widehat{C}} C_{\rho}$ --for some set of
subindexes that we call  $\widehat{C}$--, where each $C_\rho \subset C$
is a simple
$\circ$--subcoalgebra. Moreover, each $C_\rho$ admits a basis of elements
$\{t^\rho_{ij}: 1 \leq i,j \leq n_\rho\}$ such that for all $1 \leq
i,j \leq n_{\rho}$,\, $\Delta(t^\rho_{ij})= \sum_k t^\rho_{ik}\otimes
t^\rho_{kj}$
and ${t^\rho_{ij}}^\circ=t^\rho_{ji}$.
\end{enumerate}
\end{theo}
\begin{proof}

It is clear that condition \eqref{e:compactcoalg} implies condition
\eqref{e:unitinnerpro} and the proof that \eqref{e:unitinnerpro},
\eqref{e:fourierpos} and \eqref{e:prodpos} are
equivalent is the content of Observation \ref{obse:refinamiento}.

\smallskip
To prove that \eqref{e:unitinnerpro}  implies \eqref{e:decomposition}
we proceed as follows. First recall that $C$ is cosemisimple (see Observation
\ref{obse:muchascosas}) and write
$C=\bigoplus_{\rho \in \widehat{C}} C_{\rho}$ with $C_\rho$ a simple
subcoalgebra of $C$. It follows from
Corollary \ref{coro:nondeg} that each $C_\rho$ is a
$\circ$--coalgebra and from Observation
\ref{obse:basiccoeff} that the family of all $C_\rho$ for $\rho \in
\widehat{C}$ consists of all simple subcoalgebras of $C$.

Hence, we have finished the proof of first part of \eqref{e:decomposition}.
Now, call $V_\rho$ an irreducible right coideal with the property that
$\operatorname {Coeff}(V_\rho)= C_\rho$ and take a basis $\mathcal
B_\rho=\{e_1^\rho,\cdots ,e_{n_{\rho}}^\rho\}$ of
$V_\rho$ that is orthogonal with respect to the given unitary inner product
in $C$, and consider the basis $\{t^\rho_{ij}: 1 \leq i,j \leq
n_\rho\}$ of $C_\rho$ defined by the formul\ae: $\Delta(e_i^\rho)=\sum_j
e_j^\rho \otimes t_{ji}^\rho$ --see Observation \ref{obse:basic}--.
In accordance with Lemma
\ref{lema:unit_trans}, the elements
$\{t^\rho_{ij}: 1 \leq i,j \leq n_\rho\}$ satisfy the required properties:
$\Delta (t^\rho_{ij})=\sum_k t^\rho_{ik} \otimes t^\rho_{kj}$ and
${t^\rho_{ij}}^\circ=t^\rho_{ji}$.

\smallskip
Now, assuming condition \eqref{e:decomposition}
we prove condition \eqref{e:compactcoalg} as follows.
It is well known,
see for example \cite[Theorem 4.1]{kn:andres}, that being $C$ the
direct sum of simple subcoalgebras, the category of
right $C$--comodules is semisimple.
We want to prove that any $C$--comodule is
unitary and from the above it follows
that it is enough to prove that irreducible
$C$--comodules are unitary. Indeed, in order
to construct a unitary product in an
arbitrary $C$--comodule, we define the product in its irreducible
components  and then extend it by forcing
irreducible  non isomorphic components to be orthogonal.

If $V$ is an irreducible $C$--comodule, then $\operatorname{Coeff}(V)$
is of the form $C_\rho$ for some $\rho \in \widehat{C}$ and then $V
\subset C_\rho$ for some $\rho$. We
define the inner product $\langle\,\,\,,\,\,\rangle$ in $C_\rho$ setting
$\langle t^\rho_{ik},t^\rho_{j\ell}\rangle=\delta_{ij}\delta_{k\ell}$
and then by restriction we define it in $V$.
In particular it is clear that the inner product thus defined is unitary.

\end{proof}

Notice that condition (4) above, guarantees that each $C_\rho$ as a
$\circ$--coalgebra is isomorphic to $(M_n(\mathbb{C}),(-)^*)$, that is
the usual matrix coalgebra with the $\circ$--operation given by the
adjoint operator.

\begin{obse}\label{obse:cerito}
It is important to remark that the inner product constructed along the
proof of part \eqref{e:decomposition} of the above theorem, satisfies the
following additional condition:
\begin{equation} \label{eqn:cerito}\left\langle
  c^\circ,d^\circ\right\rangle =
 \overline{\langle c,d \rangle},\quad \forall c,d\in C.
\end{equation} Equivalently
 --see Observation \ref{obse:refinamiento}, (4)-- the corresponding
 associated Fourier form is symmetric.
\end{obse}

\section{The case of Hopf algebras}

\label{section:hopf}
From now on we assume that our basic $\circ$--coalgebra
$H$ has the additional structure of a Hopf algebra with a
product that is compatible with the $\circ$--structure.

\begin{defi} \label{defi:zerohopf}
Assume that $H$ is a Hopf algebra
and that $\circ$ is a $\circ$--structure on $(H,\Delta,\varepsilon)$.
We say that the pair $(H,\circ)$  is
a \emph{$\circ$--Hopf algebra} if $(xy)^{\circ}= x^{\circ}y^{\circ}$, for all $x,y \in H$.
\end{defi}

In this situation, --see \cite{kn:manin}-- the following properties are
satisfied: $1^{0}=1$, and  $S \circ S \circ = \operatorname{id}$.
In particular, it is clear that $\circ:H \rightarrow H_c$ is an
algebra automorphism and that the antipode of $H$ is an invertible
linear map  with inverse $S^{-1}=\circ S\circ$.

Next, we compare this definition with the equivalent concept of
$\star$--Hopf algebra.

Recall that a $\star$--Hopf algebra, is a Hopf algebra $H$ defined
over $\mathbb{C}$ and equipped with a conjugate linear involution $\star:H
\rightarrow H$, such that $H$:  it is a
$\star$--algebra, and the maps $\Delta:H \rightarrow H \otimes H$ and
$\varepsilon$ are $\star$--homomorphisms.

The concept of $\star$--Hopf algebra and of $\circ$--Hopf algebra are
equivalent as for a given $\star$--structure the map $S\star$ is obviously a
$\circ$--structure on $H$.

Next we collect for future reference some of the basic properties of
Fourier forms, products and integrals in the case of a Hopf
algebra.

\begin{obse}\label{obse:lem2}
The considerations that follow are consequence of
Observations  \ref{obse:refinamiento}.
If $H$ is a $\circ$--Hopf algebra and $\varphi$ is a right integral
on $H$, then
\[
\langle x,y \rangle_{\varphi}=\varphi\left( S(y^\circ)x\right),\quad \forall x,y\in H
\]
defines a unitary sesquilinear form $\langle\,\,\,,\,\,\rangle_{\varphi}:H\otimes H_c\to\mathbb{C}$.
Moreover,
\begin{enumerate}
\item
$\langle\,\,\,,\,\,\rangle_{\varphi}$ is hermitian if and only if
\(
\varphi\left( x^\circ \right)=\overline{\varphi\left( S(x)\right)},\ \forall x\in H
\),

\item
$\langle\,\,\,,\,\,\rangle_{\varphi}$ is positive definite if and only if
\(
\varphi\left( S(x^\circ)x\right)>0,\ \forall\/\, x\ne 0\in H.
\)
\end{enumerate} We prove one of the implications of part (1), the rest
of the proof is direct.
\begin{align*}
\varphi\left( S(x^\circ)y\right)
=&
\varphi \left( \left( S(x^\circ)^\circ y^\circ \right)^\circ \right)
=\\
\varphi \left( \left( S^{-1}(x) y^\circ \right)^\circ \right)
=& \overline{\varphi\left(S\left( S^{-1}(x) y^\circ \right)\right)  }
=
\overline{\varphi\left( S\left( y^\circ \right)x\right)}.
\end{align*}
\end{obse}

\begin{obse} \label{obse:dualities}
\begin{enumerate}
\item
In the situation of a Hopf algebra $H$, if $(W,\chi)\in {}^H \mathcal
M$ and
$\chi^S =\operatorname{sw}(S \otimes \operatorname{id})\chi: W \rightarrow W \otimes H$, then $(W,\chi^S)\in \mathcal M^H$.

We can consider also the standard duality  functor
$(V,\chi) \mapsto (V^*,\chi^{rS}): \mathcal M^H_f \rightarrow \mathcal M^H_f$,
that defines a duality in the category $\mathcal M^H_f$.
In explicit terms $\chi^{rS}(f)=\sum f_0\otimes
f_1$, if and only if $\sum
f_0(v) f_{1}= \sum f(v_0)S(v_1)$ where $\chi(v)=\sum v_0\otimes v_1$.

\item It follows directly from the above formula that the cannonical map
$j_V:(V,{\left(\id\otimes S^2\right)\chi})\to (V^{**},\chi^{rSrS})$
is a morphism of  $H$--comodules.

\item In the situation above, it is clear that
  $\operatorname{Coeff}\left( V^*,\chi^{rS}\right) =
  S(\operatorname{Coeff}(V,\chi))$. It is also interesting at this
  point to recall that $\operatorname{Coeff}(V^*_c,\chi^{r\circ})=
  \operatorname{Coeff}(V,\chi)^\circ$ --see Observation
  \ref{obse:coeffdual}--.
  \end{enumerate}
\end{obse}

\begin{lema}\label{lema:preparatory}
Assume that $H$ is a $\circ$--Hopf algebra and let $(V,\chi)$ be a finite
dimensional right $H$--comodule that admits an unitary inner product $\gamma :V^*\otimes V^*_c\to\mathbb{C}$. Then the map
${\psi_\gamma}_c:\left(V^*_c,\chi^{r\circ}\right)\to (V^{**},\chi^{rSrS})$
is a morphism of  $H$--comodules.
\end{lema}
\begin{proof}
In the following proof we call $\chi^{r\circ}(f)=\sum f_{(0)}\otimes f_{(1)}$ and $\chi^{rS}(f)=\sum f_{0}\otimes f_{1}$.

We want to prove that the following diagram is commutative
\begin{equation}\label{diag:cosa}
\xymatrix{
V^*_c \ar[d]_{\chi^{r\circ}} \ar[r]^-{{\psi_\gamma}_c}  & V^{**}\ar[d]^{\chi^{rSrS}} &         \\
V^*_c\otimes H \ar[r]_{{\psi_\gamma}_c\otimes\id}& V^{**}\otimes H & \\
}
\end{equation}
As $j_V:\left(V,{\left(\id\otimes S^2\right)\chi}\right)\to
\left(V^{**},\chi^{rSrS}\right)$ is a morphism --see Observation \ref{obse:dualities}-- the commutativity of \eqref{diag:cosa} is equivalent to the
commutativity of:
\begin{equation*}
\xymatrix{
V^*_c \ar[d]_{\chi^{r\circ}} \ar[r]^-{{\psi_\gamma}_c} & V^{**} \ar[r]^{j_V^{-1}} &
V\ar[d]^{\left(\id\otimes S^2\right)\chi} \\
V^*_c\otimes H \ar[r]_{{\psi_\gamma}_c\otimes\id}       & V^{**}\otimes H \ar[r]_{j_V^{-1}\otimes\id}& V\otimes H  }
\end{equation*}
Let be $f\in V^*_c$ and call $j_V^{-1}({\psi_\gamma}_c(f))=v$. Then
\begin{equation}\label{eq:fv}
\gamma(h,f)=h(v),\quad \forall h\in V^*.
\end{equation}
We have to show that
\begin{equation*}
\sum {\psi_\gamma}_c\left(f_{(0)}\right)\otimes f_{(1)}=\sum j_V(v_0)\otimes S^2(v_1).
\end{equation*}
After evaluation at an arbitrary
$g \in V^*$ and  using \eqref{eq:fv} and Observation
\ref{obse:dualities}, (1); the above equation becomes:
\begin{equation}\label{eq:2}
\sum \gamma\left(g,f_{(0)}\right) f_{(1)}=S\left(\sum \gamma(g_0,f)\,g_1\right).
\end{equation}
Let us consider a fixed $g\in V^*$ and represent it by an unique $w\in V$
such that $\overline{\gamma(g,-)}=j_V(w)$,
\textit{i.e.}:
\begin{equation}\label{eq:gw}
\gamma(g,h)=\overline{h(w)},\quad \forall h\in V^*.
\end{equation}
Hence using the characterization of $\chi^{r\circ}$ given in
Definition \ref{defi:dual}, equation \eqref{eq:2} becomes:
\begin{equation}\label{eq:3}
\sum \overline{f(w_0)}\,w_1^\circ=S\left(\sum \gamma(g_0,f)\,g_1\right)\in H.
\end{equation}
Next observe that the condition of unitarity of $\gamma$,
$\sum \gamma\left(g,f_0\right)f_1^\circ=\sum
\gamma\left(g_0,f\right)g_1$, expressed in terms of $w$ instead of $g$
becomes:
$\sum \overline{f_0(w)}\,f_1^\circ=\sum \gamma\left(g_0,f\right)g_1
$, and after applying the antipode $S$ we have:
$S\circ\left(\sum f_0(w)\,f_1\right)=S\left(\sum
  \gamma\left(g_0,f\right)g_1\right)$.
Going back to equation \eqref{eq:2} we have:

\begin{align*}
S\left(\sum
  \gamma\left(g_0,f\right)g_1\right)&= S\circ\left(\sum f_0(w)\,f_1\right)
=
\circ S^{-1}\left(\sum f_0(w)\,f_1\right) \\
&=
\circ S^{-1}\left(\sum f(w_0)\,S(w_1)\right)
=
\sum \overline{f(w_0)}\,w_1^\circ .
\end{align*}
\end{proof}

The following theorem will be of crucial importance later in the
treatment of compact involutions. Notice that this result is related
to \cite[Theorem 1.7]{kn:andrus}
and to \cite[Theorem 3.3]{kn:Larson}.
\begin{theo} \label{theo:positivity}
Assume that $H$ is a $\circ$--Hopf algebra and let $(V,\chi)$ be a finite
dimensional right $H$--comodule with the property that $V$ as well as
$V^*$ admit unitary inner products.
Then there is a positive definite isomorphism of comodules
between $(V,\chi)$ and
$\left(V,\left(\operatorname {id} \otimes S^2\right)\chi\right)$.
\end{theo}
\begin{proof}
We denote the inner products as $\beta:V \otimes V_c \rightarrow \mathbb{C}$ and
$\gamma: V^* \otimes V^{*}_c \rightarrow \mathbb{C}$ and call
$\psi_\beta: V \rightarrow V^*_c$ and $\psi_\gamma: V^*\rightarrow
V^{**}_c$ the
corresponding linear isomorphisms.

The Lemma \ref{lema:unitary} and Observation \ref{obse:dualities}
guarantee that the maps
$\psi_\beta: (V,\chi)\rightarrow (V^*_c,\chi^{r\circ})$
and $j_V:\left(V,{\left(\id\otimes S^2\right)\chi}\right)\to \left(V^{**},\chi^{rSrS}\right)$
are morphisms of $H$--comodules.
Moreover, in accordance with Observation \ref{obse:positivity},
the map $\phi$ defined by the diagram below is bijective and
positive definite with
respect to the inner product $\beta$ of $V$
\begin{equation}\label{eq:defi_fi}
\xymatrix{
V \ar[d]_{\psi_\beta} \ar[r]^-{\phi}  & V \ar[d]^{j_V} \\
V^*_c \ar[r]_-{{\psi_\gamma}_c}       & V^{**}
.}
\end{equation}

We prove that the map  $\phi:(V,\chi)\rightarrow \left(V,\left(\operatorname{id}\otimes {S}^2\right)\chi\right)$ is a morphism of $H$--comodules.

Consider the diagram that follows:
\begin{equation*}
\xymatrix{
V\ar[rrr]^{\phi}\ar[ddd]_{\chi}\ar[rd]_{\psi_\beta} & & & V\ar[ddd]^{\left(\id\otimes S^2\right)\chi}\ar[ld]_{j_V} \\
 & V^*_c \ar[d]_{\chi^{r\circ}} \ar[r]^-{{\psi_\gamma}_c}  & V^{**}\ar[d]^{\chi^{rSrS}} &         \\
 & V^*_c\otimes H \ar[r]_{{\psi_\gamma}_c\otimes\id}& V^{**}\otimes H & \\
V\otimes H \ar[rrr]_{\phi\otimes \id}\ar[ru]_{\psi_\beta\otimes\id} & & & V\otimes H \ar[lu]_{j_V\otimes\id}
}
\end{equation*}

Of the four trapezoids that appear in the diagram, the upper and lower
are commutative by \eqref{eq:defi_fi}; the left and right because
$\psi_\beta: (V,\chi)\rightarrow (V^*_c,\chi^{r\circ})$
and $j_V:\left(V,{\left(\id\otimes S^2\right)\chi}\right)\to \left(V^{**},\chi^{rSrS}\right)$ are
morphisms of $H$--comodules, respectively. As the square in the
center is commutative by the above Lemma \ref{lema:preparatory},
we conclude that the exterior rectangle is commutative, i.e. that $\phi:(V,\chi)\to \left(V,\left(\operatorname {id} \otimes S^2\right)\chi\right)$ is an $H$-comodule map.
\end{proof}

\section{Compact Quantum Groups} \label{section:compactquantum}

In this section we review some of the basic definitions related to the
concept of compact quantum group. See \cite{kn:woro} for the original
definition and also the papers cited in the Introduction for the
early development of the subject.

\begin{defi} \label{defi:basiccompactquantum} \cite[Definition 7.1, Proposition 7.5]{kn:andres} Assume
  that $(H,\Delta,\varepsilon, \mu, 1, S, \circ)$ is a $\circ$--Hopf
  algebra. We say that $H$ is a
{\em compact quantum group} if $(H,\Delta,\varepsilon,\circ)$ is
a compact $\circ$--coalgebra.
\end{defi}

In other words, a compact quantum group is a $\circ$--Hopf algebra with the
additional property that all the $H$--comodules can be endowed with a
unitary innner product.

Recall that in the situation of the equivalent set up
of $\star$--Hopf algebras --see for example
\cite{kn:koor}--, an inner product in an $H$--comodule $V$ is called
\emph{unitary} if it satisfies the following additional condition: for
all $v,w \in V$ we have that $\sum \langle v_0,w_0 \rangle
v_1w_1^{\star}= \langle v, w \rangle 1$. If we define $\circ = S \star$ this
unitary condition is equivalent to the one appearing in Definition
\ref{defi:unitary}. Now, the definition of compact quantum group in
terms of a $\star$ operation is that of a $\star$--Hopf algebra with
the property that all the $H$--comodules admit a \emph{unitary} inner
product. Hence the definition of compact quantum group
apperaring for example in
\cite{kn:andrus0}, \cite{kn:andrus} or \cite{kn:koor} in terms of a
$\star$--structure on $H$,  is equivalent to the one we presented above.

As we mentioned before, the above Definition \ref{defi:basiccompactquantum}
of compact quantum group --or the equivalent concept in
terms of a $\star$--operation-- is the algebraic counterpart of the
original definition by Woronowicz, that is obtained from the one above
by a process of ``completion''. We refer the reader to
\cite[Proposition 1.4]{kn:andrus0},
\cite[Chapter 2, Section 4]{kn:guichardet} or \cite{kn:koor2} for the
precise description of this
relationship.

We recall some results concerning cosemisimple Hopf algebras and then
look in the case of $\circ$--Hopf algebras at the behavior
of the integral in relation to the $\circ$--operator.

\begin{obse} \label{obse:int1}
\begin{enumerate}
\item Suppose that $H$ is a cosemisimple Hopf algebra.
We can consider $\mathbb{C}$ as a right --or left--
$H$--comodule with the trivial structure given by the unit. In
  Observation \ref{obse:basiccoeff} we noticed that we can write
  $H=\bigoplus_{\rho \in \widehat{H}}H_\rho$ where $H_\rho$ is the
family of all simple subcoalgebras of $H$. In this situation we call
$\rho_0$ the element in $\widehat H$, with the property that
$H_{\rho_0}=\mathbb{C}$. Let $\varphi:H \rightarrow \mathbb{C}$ be
the projection from $H$ onto the component
  $H_{\rho_0}$. It is clear that $\varphi$ is a morphism of right
  and left comodules and this implies that
it is a right and left normal integral on $H$.
\item In particular, the antipode $S$ permutes the simple
  subcoalgebras $H_\rho$ for $\rho \in \widehat H$ and leaves $H_{\rho_0}=
  \mathbb{C}$ fixed.

\item In this situation if $\varphi: H \rightarrow \mathbb{C}$
is the normal --right and left-- integral on
$H$, we can consider in $H$ the Fourier product $\star_\varphi$ and
the normal Fourier form $\omega_\varphi$
induced by $\varphi$ --see Observation \ref{obse:lem2}--.

\[
x \star_\varphi y=\sum \varphi(S(x)y_1)y_2=\sum x_1\varphi(S(x_2)y),\quad \forall x,y\in H.
\]
\[
\omega_\varphi(x,y)=\varphi(S(x)y).\]

The Fourier form considered above is left and right non degenerate and
satisfies:
\[\omega_\varphi\left(x , y\right)=
\omega_\varphi\left(y, S^2(x)\right), \quad \forall x,y\in H_\rho.\]

See for example \cite{kn:sweedler2}, for the non degeneracy of the form.

\item Moreover, in the case of a $\circ$--Hopf algebra,
the map $\circ:H \rightarrow
  H$ leaves invariant all the $H_\rho$, including
  $H_{\rho_{0}}$--see Lemma \ref{lema:unitary}--.
Hence it follows that $\varphi\/\, \circ = \overline
  {\varphi}$, \textit{i.e.}, $\varphi(x^\circ)=\overline{\varphi(x)}$
  for all $x \in H$.

\item The bilinear map $\langle\,\,\,,\,\,\rangle_\varphi$ ,$\langle
  x,y \rangle_\varphi=\varphi(S(y^{\circ})x)$, for all $x,y \in H$, is
  hermitian, non degenerate and unitary. The hermitianity is equivalent to the
  property that $\varphi(x^\circ)= \overline {\varphi(S(x))}$ --see
  Observation \ref{obse:lem2}-- that is a
  consequence of the properties (1) and (4) above.

\end{enumerate}
\end{obse}

\begin{obse} \label{obse:equalities2}
\begin{enumerate}
\item It follows from Theorem
  \ref{theo:positivity} that for an arbitrary simple coalgebra $C
  \subset H$ of the compact quantum group $H$,
  we have that $S^2(C)=C$.
Indeed, in accordance with the mentioned theorem, if $V$ is a finite
dimensional $H$--comodule, there is an isomorphism of
$H$--comodules
$\varphi: (V,\chi) \rightarrow (V,(\operatorname {id} \otimes
S^2)\chi)$. Hence, $\operatorname {Coeff}(V)
= S^2(\operatorname {Coeff}(V))$ and taking $V$
an irreducible comodule that has $C$ as coalgebra of
coefficients, we deduce that $C= S^2(C)$. Notice that this
result is valid for any simple subcoalgebra $C$ of a Hopf
algebra --not necessarily with a $\circ$--structure--
and was proved by Larson in
\cite{kn:Larson}[Thm. 3.3].
\item For future reference we write down the following formula, valid
  in our context for all $x,y \in H$
\begin{equation}\label{eqn:fourier2}
\sum x_1 \langle x_2,y \rangle_\varphi = \sum
  S^2(y_1^\circ) \langle x,y_2 \rangle_\varphi
\end{equation}

Explicitly, we need to show that: $\sum x_1 \varphi(S(y^\circ) x_2)=
\sum S^2(y_1^\circ) \varphi(S(y_2^\circ) x)$, that by a
change of variables becomes: $\sum x_1 \varphi(y x_2)= \sum S(y_1)
\varphi(y_2 x)$. This last equality can be proved directly  using
the fact that $\varphi$ is an integral --see \cite{kn:ferrerfourier}--.
\end{enumerate}
\end{obse}

Let $H$ be a compact quantum group, next we
prove that the non degenerate unitary
form associated to the integral --that is denoted as
$\langle\,\,\,,\,\,\rangle_\varphi$-- and  given explicitly as $\langle
x,\,y\rangle= \varphi(S(y^\circ)x)$ is in fact an inner
product. This is the content of Theorem \ref{theo:positivity2}.

In order to prove this positivity result we will use
  the {\em unitary inner product}
whose existence was established along the proof
of Theorem \ref{theo:coalgcompacta}.  This product will be denoted as
$\langle\,\,\,,\,\,\rangle$ and it satisfies the equation $\langle
x^\circ,\,y^\circ\rangle=\overline{\langle x,y \rangle}$ --see
Observation \ref{obse:cerito}, \eqref{eqn:cerito}--.

As the Fourier forms associated to
$\langle\,\,\,,\,\,\rangle_\varphi$ and
$\langle\,\,\,,\,\,\rangle$ are unitary, the subcoalgebras $H_\rho$ are
orthogonal with respect to both. The mentioned positivity result in $H$
will be deduced of the corresponding result for each
$H_\rho$.  Explicitly the Fourier products
$\star$ and $\star_\varphi$ associated respectively with
$\langle\,\,\,,\,\,\rangle$ and  $\langle\,\,\,,\,\,\rangle_\varphi$
are given as:
\begin{align*}
x\star y &= \sum \left\langle y_1,x^\circ \right\rangle y_2 = \sum x_1\left\langle y,x_2^\circ \right\rangle, \\
x\star_\varphi y &= \sum \left\langle y_1,x^\circ \right\rangle_\varphi y_2
= \sum x_1\left\langle y,x_2^\circ \right\rangle_\varphi
= \sum\varphi(S(x)y_1)y_2=\sum x_1\varphi(S(x_2)y).
\end{align*}
Next Lemmas \ref{lema:lem3}, \ref{lema:lem4} are preparatory for the
proof ot the positivity result mentioned above.

\begin{lema}\label{lema:lem3} Let $H$ be a compact quantum group and
  let $H_\rho$ be one of its simple components as before.
If $H_\rho=\bigoplus_{a} R_a$ is the decomposition of $H_\rho$ into
irreducible right coideals (that are pairwise
$\langle\,\,\,,\,\,\rangle$--orthogonal in accordance with Corollary
\ref{coro:nonisoorth}),
then:
\begin{enumerate}
\item $H_\rho=\bigoplus_a (R_a)^\circ$ is also a decomposition of
  $H_\rho$ into pairwise orthogonal subspaces.
\item The element $s_\rho \in H_\rho$, neutral for the restriction of
$\star$ to $H_\rho$, is $\circ$--stable, i.e., $s_\rho=s_\rho^\circ$.
\item If we write $s_\rho=\sum_a s_a$,
with $s_a \in R_a$ then, for all $x \in R_a$, $x=s_a\star x$ and if $b \neq a$, then $s_b\star
  x=0$. In particular $s_a\ne 0$, for all $a$.
\end{enumerate}
\end{lema}
\begin{proof}
(1) The properties of this decomposition can be easily deduced from
the fact that
$R_a$ is orthogonal to $R_b$ when $a\neq b$ --see Corollary
\ref{coro:nonisoorth}--(we are also using that for $x,y \in H_\rho$,
$\langle x^\circ\,,y^\circ\rangle=\overline{\langle x\,,y\rangle}$,
c.f. Observation \ref{obse:cerito}).

(2) The existence of $s_\rho$ follows from Observation
\ref{obse:lem1}. We have that:
$s_\rho^\circ=s_\rho^\circ \star s_\rho= \left (s_\rho^\circ \star
  s_\rho \right)^\circ= \left (s_\rho^\circ \right)^\circ=s_\rho$
--see Observation \ref{obse:refinamiento}--.

(3)
From the equality: $x=s_\rho\star x= \sum_b (s_b \star x)$ and the fact that $s_b\star x
\in R_b$ for all $b$ --see Observation \ref{obse:lem1}, we deduce our result.
\end{proof}

\begin{lema}\label{lema:lem4}
Let $H$ be a compact quantum group,
  let $H_\rho$ be one of its simple components and let
$H_\rho=\bigoplus_{a} R_a$ be the decomposition of the simple
subcoalgebra $H_\rho$ into irreducible right coideals.
Then, for each $a$ we call
$A_a:R_a \rightarrow R_a$ the linear map defined by the equality:
\[
A_a(x)=s_a \star_\varphi S^2(x), \quad \forall x\in R_a.
\]
Then:
\begin{enumerate}
\item
For all $x,y \in R_a$, we have
\begin{equation}\label{claim1}
\left\langle y^\circ,x^\circ \right\rangle_\varphi= \langle A_a(x),y\rangle.
\end{equation}
\item
The map $A_a$ is a positive operator with respect to the inner product
$\langle\,\,\,,\,\,\rangle$.
\end{enumerate}
\end{lema}
\begin{proof}
(1)
To prove the equality \eqref{claim1}, we take $x,y \in R_a$ and
deduce that:
\begin{align*}
\left\langle y^\circ,x^\circ \right\rangle_\varphi
&=
\varepsilon\left( x \star_\varphi y^\circ \right)
=
\varepsilon\left( y^\circ \star_\varphi S^2(x) \right)
=
\varepsilon\left( y^\circ \star s_\rho \star_\varphi S^2(x) \right)
=
\left\langle s_\rho \star_\varphi S^2(x),y \right\rangle \\
&=
\sum_b\langle s_b \star_\varphi S^2(x),y\rangle
=
\langle A_a(x),y\rangle .
\end{align*}

The second in the chain of equalities above comes from the fact that:
\[\varepsilon\left( x \star_\varphi
  y^\circ \right)= \omega_\varphi(x\,,y^{\circ})=  \omega_\varphi(y^\circ\,,S^2(x))=
\varepsilon\left( y^\circ \star_\varphi S^2(x) \right),\] --see Observation \ref{obse:int1} part (3)--.

(2)
The map $A_a:(R_a, (\id\otimes S^2)\Delta)\to (R_a,\Delta)$ is an $H$-comodule map. Indeed,
\begin{align*}
\sum A_a(x_1)\otimes S^2(x_2) &=
\sum s_a \star_\varphi S^2 (x_1)\otimes S^2(x_2)=
\sum s_a \star_\varphi S^2 (x)_1\otimes S^2(x)_2 =\\
&(\small{{\text{see Def. \ref{defi:transform}}}}) =
\Delta \left(s_a \star_\varphi S^2 (x) \right)
=\Delta(A_a(x)).
\end{align*}
From Theorem \ref{theo:positivity} and Schur's Lemma we deduce that there is a complex scalar $\gamma_a$ and
a $\langle\,\,\,,\,\,\rangle$--positive operator $P_a$ on $R_a$ with the property that $A_a=\gamma_a P_a$.
Let $B_a=\{e_1,e_2,\cdots,e_{n}\}$ be a $\langle\,\,\,,\,\,\rangle$--orthonormal basis of $R_a$.
If we write $\Delta(s_a)=\sum _i e_i \otimes h_i$, then
\begin{align*}
\sum_i e_i \langle e_j, h_i^\circ \rangle  & =  s_a \star e_j
= e_j.
\end{align*}
so we deduce $\langle e_j,h_i^\circ \rangle= \delta_{i,j}$ and
therefore $h_i=e_i^\circ$. Then,
$\Delta(s_a)=\sum_i e_i \otimes e_i^\circ $ and --see Observation
\ref{obse:lem1}--
$
\sum_i e_i \star_\varphi e_i^\circ = s_a.
$
Applying $\varepsilon$ to this last equality we deduce that:
$\sum_i \langle e_i^\circ, e_i ^\circ
\rangle_\varphi=\varepsilon(s_a)$.

The positivity of $\gamma_a$ follows from the following computation:
\begin{align*}
\gamma_a \sum_i \langle P_a(e_i),e_i\rangle = & \sum_i \langle
A_a(e_i),e_i\rangle = \sum_i \langle e_i^\circ, e_i ^\circ \rangle_\varphi=\varepsilon(s_a)
= \\ = \varepsilon(s_\rho \star s_a)= & \langle
s_a,s_\rho^\circ\rangle = \langle s_a, s_\rho\rangle= \langle s_a, s_a\rangle >0.
\end{align*}
Hence $\gamma_a$ is a strictly positive real number and thus
$A_a$ is a positive operator.
\end{proof}

\begin{obse}\label{obse:obse4}
In particular for all $x \in R_a^\circ$,
\begin{equation}\label{claim2}
\left\langle x,x \right\rangle_\varphi = \langle A_a(x^\circ),x^\circ\rangle.
\end{equation}
\end{obse}

\begin{theo}\label{theo:positivity2}
If $H$ is a compact quantum group and $\varphi: H \rightarrow
\mathbb{C}$ is the associated normal integral, then the hermitian form $\langle\,\,\,,\,\,\rangle_\varphi$ is positive definite.
\end{theo}
\begin{proof}
As we proved in Theorem \ref{theo:coalgcompacta} there is a
decomposition $H=\bigoplus_{\rho\in \widehat H}H_{\rho}$, where
$H_\rho$ are non isomorphic simple
subcoalgebras that are also $\circ$--stable and $\langle\,\,\,,\,\,\rangle_\varphi$--orthogonal.
For each $\rho\in \widehat H$ we have a decomposition
$H_\rho=\bigoplus_{a} R_a$, into irreducible right
coideals that are all isomorphic --see Lemma \ref{lema:lem3}--.
Then, $H_\rho =\bigoplus_{a} R_a^\circ$ and we claim that the left coideals
$R_a^\circ$ are $\langle\,\,\,,\,\,\rangle_\varphi$--orthogonal.

\medskip

As $R_a$ and $R_b^\circ$ are respectively right and left
coideals, from Observation \ref{obse:lem1} we deduce that
$R_a\star R_b^\circ\subset R_a\cap R_b^\circ$.
Conversely, if $x\in R_a \cap R_b^\circ$, then $x=s_a\star x\in R_a\star
R_b^\circ$ and we deduce that $R_a\star R_b^\circ = R_a\cap R_b^\circ$.

Applying $\varepsilon$ the the above equality we obtain:
\[
\varepsilon\left(R_a\cap R_b^\circ\right) =
\varepsilon\left(R_a \star R_b^\circ\right) =
\left\langle R_b^\circ , R_a^\circ\right\rangle =
\overline{\left\langle R_b , R_a \right\rangle} = 0,
\] where we used above that for all $x,y \in H$:
$\langle x^\circ, y^\circ \rangle = \overline{\langle x, y \rangle}$ --see
Observation \ref{obse:cerito}, \eqref{eqn:cerito}--.
Now, if $x\in R_b^\circ$ and $y\in R_a^\circ$, then
$y^\circ\star_\varphi x\in R_a\cap R_b^\circ$ and
$\langle x,y\rangle_\varphi=\varepsilon(y^\circ\star_\varphi x) \in
\varepsilon(R_a \cap R_b^\circ)=\{0\}$.

Once that the above $\langle\,\,\,,\,\,\rangle_\varphi$--orthogonality
is established, it is enough to show that
$\langle x,x\rangle_\varphi >0$, for all $a$ and $\forall x \in
R_a^\circ,\ x\neq 0$,  and this is exactly the content of
Lemma \ref{lema:lem4} --see also Observation \ref{obse:obse4}--.
\end{proof}

A proof in matricial terms of the above theorem, appears in
 \cite[Theorem 2.8]{kn:koor} --see also the previous paper
 \cite{kn:andrus0}--.

\begin{coro} \label{coro:positiveS2}
If $H$ is a compact quantum group, then the
automorphism $S^2:H \rightarrow H$ is positive definite with respect
to the inner product defined by the integral $\varphi$.
\end{coro}
\begin{proof} The result follows from the equalities:
$\langle S^2(x),x \rangle_\varphi=\varphi(S(x^\circ)S^2(x))=\varphi(S(x) x^\circ)=\langle
  x^\circ,x^\circ\rangle_\varphi$.
\end{proof}

We summarize the above results in the following theorem (see
\cite[Proposition 2.4]{kn:andrus}).
\begin{theo} \label{theo:maincarac} A
$\circ$--Hopf algebra $H$ is a compact quantum group if and only if it is
  cosemisimple and the normal integral $\varphi:H \rightarrow
  \mathbb{C}$ is positive definite in the sense
that $\varphi(S(x^{\circ})x)
  >0$ for all $0 \neq x \in H$.
\end{theo}

Notice that in accordance with Observation \ref{obse:int1} part (5), the
associated bilinear map $\langle\,,\,\rangle_\varphi$ is always
hermitian non--degenerate and unitary.
\section{The conjugacy of compact involutions}
\label{section:compactinv}
In this section we present a detailed proof --using the methods introduced
above-- of a result that appeared in \cite[Theorem 2.6]{kn:andrus}:
if a cosemisimple Hopf algebra has two different compact
$\circ$--structures,
then they are conjugate by a positive Hopf algebra
automorphism.

The proof is self--contained and should be thought as an
expansion and elaboration of the one appearing in \cite{kn:andrus}. In
the same manner than
the original
proof, ours is based in a few simple facts
about positive linear
transformations on finite dimensional complex vector spaces --and
locally finite algebras-- equipped with an
inner product.

We recall some of the needed elementary results below.

\begin{obse} \label{obse:linearalgebra} Let $V$ be a
vector space equipped
  with an inner product $[\,\,\,,\,\,]$ and $N:V \rightarrow V$ be a
  positive linear transformation. Assume moreover that $V$ is
  $N$--locally finite in the sense that $V$ is the direct sum of
  $N$--stable finite dimensional spaces.
It is well known that if we work on the finite dimensional $N$--stable
pieces of $V$, we can find
positive real numbers $\lambda_i:i=1,\cdots,k$ and orthogonal
projections $E_i:i=1,\cdots,k$ such that
$N=\lambda_1 E_1 \oplus \cdots \oplus \lambda_k E_k$.
The operator $P=(\lambda_1)^{1/2}
E_1 \oplus \cdots \oplus (\lambda_k)^{1/2}E_k$ is the only
positive operator with respect to $[\,\,\,,\,\,]$ satisfying that
$P^2=N$. The inverses of $N$ and $P$ are also positive operators with
respect to $[\,\,\,,\,\,]$. All this can be globalized to all of $V$
in the obvious manner.

Moreover, if there exists a linear involution $\sigma:V \rightarrow V$
such that $\sigma N = N^{-1}\sigma$, then
  $\sigma P = P^{-1}\sigma$, and this can be proved locally.
In that situation, if $\lambda$ is an eigenvalue of
  $N$ and $E_\lambda$ is the corresponding eigenspace, then from the
  equality $\sigma N \sigma = N^{-1}$ we deduce that
  $\sigma(E_\lambda)=E_{\lambda^{-1}}$. Then, $P\sigma|_{E_\lambda}=
  (\sqrt{\lambda})^{-1}\sigma|_{E_\lambda}$ and thus, for
$y \in E_\lambda$, $\sigma P \sigma y =
  (\sqrt{\lambda})^{-1} y= P^{-1}y$.
\end{obse}

\begin{obse} \label{obse:multiplicativity}
Here --for later use--
we present two applications of the above construction, one to
algebras and the other to simple coalgebras.

\begin{enumerate}
\item Assume that additionally to the situation above we have that  $V=A$ is a $\mathbb{C}$--algebra equipped with an inner
product $[\,\,\,,\,\,]$ and
that $N:A \rightarrow A$ is a positive
algebra automorphism.
Using the fact that $N$ and $P$ have the same eigenspaces,
it is easy to prove that $P$ is an algebra automorphism. Clearly, from
$N(1)=1$ we deduce that $1 \in H$ is an eigenvector and then
$P(1)=1$. Moreover, it is easy --and enough--
to check the multiplicativity for a pair $x,y $ of eigenvectors of $P$
--or equivalently, eigenvectors of $N$--. As $xy$ is also an
eigenvector of $N$ --and then of $P$--, the multiplicativity follows
immediately.
\item Later we will use also a dual version of the above result.
If $N$ is a positive multiplicative coalgebra
  automorphism in the same situation that above, then its positive
  square root is also a coalgebra automorphism.
\item If $C$ is an arbitrary coalgebra, it
  is clear that a linear map $W: C \rightarrow C$ is a morphism of
  right comodules if and only if there exists a functional $\omega: C
  \rightarrow \mathbb C$ with the property that $W=w \star
  \operatorname{id}$. In that case $\omega=\varepsilon W$.
Assume moreover that $C$ is a simple coalgebra and that
$W$ is positive with respect to
a given inner product $\langle
\,\,\,, \,\,\rangle$ on $C$. We want to prove that in this situation
there exists another
functional $\rho:C \rightarrow \mathbb C$ such that:
$\rho^2=\omega$ and the map $\rho \star \id: C \rightarrow C$ is
positive.
Indeed,  if we write the decomposition of $C$ as a sum
of $W$--eigenspaces: $C=\bigoplus
  E_a$, then $W|_{E_a}=\lambda_a
  \operatorname{id}$ with $\lambda_a$ a positive real number. In that
  case $\omega|_{E_a}= \lambda_a \varepsilon$. If we
  call $P$ the positive square root of $W$, it is clear that on each
  $E_a$, the operator $P$ is of the form $\sqrt{\lambda_a}
  \operatorname {id}$. From this expression it follows that
$P$ is also a morphism of right comodules.
To prove this we take $x \in E_a$ and hence,
  $P(x)=\sqrt{\lambda_a}x$ and $W(x)=\lambda_a x$. As $W$ is a morphism
  of right comodules, we deduce that if we write $\Delta(x)=\sum u_i
  \otimes v_i$ with vectors $v_i$ linearly independent,
then $\lambda_a \sum u_i
  \otimes v_i = \sum W(u_i) \otimes v_i$. Hence, $W(u_i)= \lambda_a
  u_i$ --i.e. $\Delta(E_a) \subset E_a \otimes C$--
and by the definition of $P$,
we know that on the eigenvectors $u_i$, $P(u_i)=
  \sqrt{\lambda_a}u_i$. Hence $P$ is a morphism of right comodules and
  as such $P=\rho \star \operatorname {id}$ for a certain linear
  functional $\rho: C \rightarrow \mathbb C$. In this situation it is
  clear that $\rho^2=\rho \star \rho = \omega$.
\item Observe that we have
  proved above in part (3), that if $V$ is a right $C$--comodule, $W: V \rightarrow V$
  is a morphism of $C$--comodules and $E$ is an eigenspace of
  $W$, then $E \subset V$ is a $C$--subcomodule. There is also a
  version of this result for left $C$--comodules.
\end{enumerate}
\end{obse}

\begin{theo}\label{theo:conmutation}
Let $(H,\circ)$ be a compact quantum group and $\lozenge:H \rightarrow H_c$ be an involution
such that $(H,\lozenge)$ is a $\circ$--Hopf algebra. Then:
\begin{enumerate}
\item There exists a Hopf algebra automorphism $P$, positive with
respect to the inner product associated to $\circ$, such that the
involutions $\lozenge$ and $P \circ P^{-1}$ commute with each
other.
\item If $(H,\lozenge)$ is also a compact quantum group,
then $\lozenge = P \circ P^{-1}$.
\end{enumerate}
\end{theo}
\begin{proof}
\begin{enumerate}
\item Consider the map $Q:H \rightarrow H$ defined by $Q=\lozenge
\circ$. Clearly, $Q^{-1}= \circ \lozenge$, $Q$ is an automorphism of
Hopf algebras --therefore it commutes with the antipode-- and $\varphi
Q= \varphi$ where $\varphi$ is the normal
integral on $H$. Moreover, $Q \circ =\circ Q^{-1}$ and $\lozenge Q=
Q^{-1} \lozenge$.
First, we prove that $Q$ is selfadjoint with respect to the inner product
associated to $\circ$. Indeed:
\begin{align*}
&\langle Q(x), y \rangle = \varphi(S(y^\circ) x^{\circ
\lozenge}) = (\varphi Q^{-1})(S(y^\circ) x^{\circ
\lozenge})\\
= &\varphi(S(y^{\circ \lozenge \circ}) x) = \varphi
(S(Q(y)^\circ) x) = \langle x, Q(y)\rangle .
\end{align*}
 Therefore, $Q^2$ is a positive
automorphism of Hopf algebras and also $Q^2 \,\circ = \circ\, Q^{-2}$
and $\lozenge\, Q^2= Q^{-2}\, \lozenge$.

Let $H'_\rho=
\left \{
\begin{array}{ll}
H_\rho &\mbox { if } H_\rho^\lozenge=H_\rho\\
H_\rho \bigoplus H_\rho^\lozenge &\mbox{ if not}
\end{array}
\right .
$\\
Clearly, for each $\rho \in \widehat {H}$, $H'_\rho$ is a $\lozenge$-subcoalgebra of $H$ and
$H=\bigoplus_{\rho \in \overline H} H'_\rho$, for some $\overline H
\subseteq \widehat {H}$.

Then, $H'_\rho$ is $Q^2$-invariant and finite
dimensional and for each $\rho \in \overline H$ we can define
the operator $P_\rho: H'_\rho \rightarrow H'_\rho$, that is the
positive fourth-root of $Q^2|_{H'_\rho}$.
We define the operator $P:H \rightarrow H$
as the sum of all these $P_\rho$.
In this situation $P$ is invertible, positive and
$P^4=Q^2$.

Now, by applying twice the Observation \ref{obse:linearalgebra} we
deduce that $P\, \circ = \circ\, P^{-1}$ and $\lozenge\, P= P^{-1}\,
\lozenge$ .

Then:
\begin{align*}
\lozenge\, P \circ P^{-1}\, \lozenge\, P \circ P^{-1} =&\, \lozenge \circ
P^{-2}\, \lozenge\, P \circ P ^{-1} =\\ \lozenge \circ \lozenge\, P^3
\circ P^{-1}=&\, Q^2P^{-4}=\operatorname{id}_H.
\end{align*}
Hence, we conclude that
$\lozenge\, P \circ P^{-1}=(\lozenge\, P \circ\, P^{-1})^{-1}= P\,
\circ\, P ^{-1}\, \lozenge$.

\item
Assume now that $(H,\lozenge)$ is also a compact quantum group and call $[\,\,\,,\,\,]$ the inner product associated to $\lozenge$, i.e. defined by
  $[x,y]=\varphi(S(y^\lozenge)x)$.\\
In this situation, if $x,y \in H$, we have that $\langle Q^{-1}x,y\rangle=
  \varphi(S(y^\circ)Q^{-1}x)= \varphi(S^{-1}(y)^{\circ} Q^{-1}x)=
(\varphi Q)(S^{-1}(y)^{\circ} Q^{-1}x)=
\varphi(S^{-1}(y)^\lozenge x)=\varphi(S(y^\lozenge)x)= [ x,y]$. Hence $Q^{-1}$ is positive on $(H,\langle\,\,\,,\,\,\rangle)$ and so is $Q$.\\
Now, as both $P^2$ and $Q$ are positive square roots of $P^4=Q^2$, we
deduce that  $P^2=Q$ and:
$$
P \circ P^{-1}= \circ\, P^{-2} = \lozenge \lozenge \circ P^{-2} = \lozenge Q P^{-2}=\lozenge
$$
and the proof is finished.
\end{enumerate}
\end{proof}
\begin{obse} \label{obse:existence} The general problem of the
  \emph{existence} of a compact involution for a cosemisimple Hopf
  algebra is --in the knowledge of
  the authors-- wide open.

Due to well known results on
  semisimple Hopf algebras --see for example \cite{kn:schneider}--,
in the case that the original Hopf
  algebra $H$ is finite dimensional, a compact finite quantum
  group is simply a (finite dimensional) semisimple Hopf algebra
(with a normal left and right integral that we call $\varphi:H
\rightarrow \mathbb{C}$) with
  the property that the sesquilinear form in $H$ defined by the formula:
$\langle x,y \rangle = \varphi(S(y^{\circ}) x)$, is positive
definite.

In the survey article
\cite{kn:andrus1}, the following Question 7.8 is raised.

\emph{Given a semisimple Hopf algebra H, does it admit a compact
  involution?}.

As far as the authors are aware of, the answer to this
question is not known.

Regarding this point, the
particular case of abelian extensions of Hopf algebras is
considered in \cite{kn:masuoka}. In particular when the original Hopf algebras are of the form $H=\mathbb{C}F$ and $K=\mathbb{C}^G$, the author briefly
considers the extension problem for $H$ and $K$ considered as {\em
  finite quantum groups} endowed with their natural $\circ$--structures.
He proves that the group of extensions of $H$ and $K$ as $\circ$--Hopf
algebras coincides with the group of extensions as Hopf algebras
(\cite[Remark 1.11, (2)]{kn:masuoka}).

One expects that this observation could be generalized to the case of
arbitrary finite dimensional (abelian) extensions.

In the case of a not necessarily finite dimensional (cosemisimple) Hopf
algebra $H$, the classical situation could be an indication of the
results to be expected. At most one should expect to
construct a compact involution in a sufficiently large {\em quotient}
of $H$.

\end{obse}
\section{The mien of the antipode of a compact quantum group}
\label{section:mien}

In this section we intend to obtain a better understanding of the
antipode antiautomorphism $S$ in the case of a compact quantum group,
and for this purpose
we apply in this particular set up,
some of the standard tools in the theory of coFrobenius Hopf algebras:
the modular function, the Nakayama automorphism, Radford's formula,
etc. Our algebraic methods
differ from the more standard analytical methods applied for example
in \cite{kn:mnw}.

\medskip

We start by defining the positive antipode.
Once we know (see Corollary \ref{coro:positiveS2}) that
$S^2$ is a positive Hopf algebra automorphism with
respect to the inner product $\langle\,\,\,, \,\,\rangle$ given by the
normal integral, it follows from the results mentioned in
Observation \ref{obse:multiplicativity}  that there exists a positive
automorphism of algebras and coalgebras that is the positive square
root of $S^2$.

\begin{defi} \label{defi:posantip}
Let $H$ be a compact quantum group with integral $\varphi$ and
associated inner product $\langle\,\,\,,\,\,\,\rangle$.
Then the unique positive automorphism of
algebras and coalgebras $S_+:H \rightarrow H$  that satisfies
\begin{equation}\label{eqn:posantipode}
S^2=S^2_+.
\end{equation}
is called the
{\em positive antipode} of $H$.
\end{defi}

\begin{obse}
It is clear that  $S_+$ leaves all the
simple components $H_\rho$ invariant.
\end{obse}

Our next objective is to prove that $S_+$ is an inner automorphism and
then to deduce that $S^2: H \rightarrow H$ is also an inner
automorphism --see \cite{kn:mnw} for a proof of this result--.

This result will follow easily from the following well known fact:
a coalgebra
automorhism of a matrix coalgebra is inner
--the dual of this result is attributed to A.A. Albert in a more
general version--.

\begin{theo}[A.A.Albert]\label{theo:albert} Let $V$ be a finite dimensional vector
space and consider the coalgebra $c(V)$. Let $T: c(V) \rightarrow
c(V)$ be an automorphism of coalgebras. Then, there exists a
convolution invertible functional
   $\tau:c(V) \rightarrow \mathbb{C}$, such that $T=(\tau \star
   \operatorname {id}) (\operatorname
   {id} \star\, \tau^{-1})$. Moreover, the linear functional $\tau$ is unique up to
   multiplication by a non zero scalar.
 \end{theo}

\begin{obse}\label{obse:inverse}
Assume that $C$ is a
  finite--dimensional coalgebra and that $T:C \rightarrow C$ is a linear map
  that can be written as: $T=(\tau \star
  \operatorname{id})(\operatorname{id} \star\, \rho)$ and that satisfies the
  equality $\varepsilon T = \varepsilon$. Then $\rho=\tau^{-1}$,
  $\tau T=\tau$, $\rho T=\rho$ and $T$ is a morphism of coalgebras.
\end{obse}

The following consequences of Albert's theorem \ref{theo:albert},
will be used later.
\begin{coro}\label{coro:albertpos} Assume that $C$ is a
  $\circ$--coalgebra that is simple and equipped with a unitary inner
product $\langle\,\,\,,\,\,\,\rangle$. Let
$T$ be an automorphism of coalgebras and suppose that it can be
  decomposed as $T = (\tau \star
  \operatorname {id})(\operatorname {id} \star\, \tau^{-1})$ where
$\tau \star \operatorname {id}$ as well as
$\operatorname {id} \star\, \tau^{-1}$ are selfadjoint.
If $T$ is positive then $\tau \star
  \operatorname {id}$ and
$\operatorname {id} \star\, \tau^{-1}$ can be taken to be positive.
\end{coro}
\begin{proof}
Notice that the existence of a functional $\tau$ yielding a
decomposition as above follows easily from Theorem
\ref{theo:albert}. The selfadjointness has to be assumed as hypothesis.

Let us call
$E_\lambda$ and $F_\mu$ the  eigenspaces for $\tau \star\,
\operatorname{id}$ and $\operatorname {id} \star\,
\tau^{-1}$ corresponding to the eigenvalues $\lambda$ and
$\mu$ respectively.

As $E= E_\lambda \cap F_\mu$ is a common eigenspace for the
--commuting-- operators
$\tau \star\, \operatorname {id}$ and
$\operatorname {id} \star\,\, \tau^{-1}$, it is clear that $\lambda \mu$
is an eigenvalue of $T$ and then it is positive and we deduce that
$\lambda$ and $\mu$ --are real and-- have the same sign.

Next we prove that one of the
subspaces
$E_{+}=\bigoplus_{\lambda >0 ,\,\mu > 0} E_{\lambda} \cap F_\mu$ or
$E_-=\bigoplus_{\lambda <0,\,\mu <0} E_{\lambda} \cap F_\mu$
has to be trivial. Observe that they
satisfy that $C=E_+ \oplus E_-$.

It follows from Observation \ref{obse:multiplicativity} part
(4) that $\Delta(E_\lambda) \subset E_\lambda
\otimes C$ and similarly $\Delta(F_\mu) \subset C \otimes
F_\mu$. Then $\Delta(E_\lambda \cap F_\mu) \subset (E_\lambda
\otimes C) \cap (C \otimes F_\mu) = E_\lambda \otimes F_\mu$. It is
clear that for $\lambda >0,\,\mu >0$, $E_\lambda\,\text{and}\, F_\mu \subset E_+$
--and similarly for negative eigenvalues-- and then it follows
that: $\Delta(E_+) \subset E_+ \otimes E_+,\,\,\Delta(E_-)
\subset E_- \otimes E_-$. As the coalgebra $C$ is simple, we deduce that
one of these subspaces (subcoalgebras) has to be trivial.

Hence, after eventually changing $\tau$ by $-\tau$ we deduce that we
can assume that the functional appearing in the formula $T = (\tau \star
  \operatorname {id})(\operatorname {id} \star\, \tau^{-1})$ has the
  property that $\tau \star
  \operatorname {id}$ and $\operatorname {id} \star\, \tau^{-1}$ are
  positive.

\end{proof}

\begin{lema}\label{lema:adjoint} Assume that $H$ is a compact quantum
   group and that $\langle\,\,\,,\,\,\rangle$ is the inner product associated
   to the normal integral $\varphi$.
 For $\tau \in H^*$, the adjoints of the maps $\id \star \tau \,,\,
 \tau \star \id: H \rightarrow H$ with
   respect to the inner product are:
 $(\operatorname {id} \star\, \tau)^*= \id \star
   (\overline{\tau}\, \circ)$ and $(\tau \star \id)^*= (\overline{\tau}S^2\circ)
   \star\,\, \id$ respectively.
 \end{lema}
 \begin{proof} We only prove the second assertion, the first is similar. We
   take $z,w \in H$ and compute explicitly $\langle (\tau \star \id)(z),
   w\rangle$ and $\langle z, ((\overline{\tau}S^2\circ)\star \id) (w)\rangle$.
 In this situation the first expression is : $\langle (\tau \star \id)(z), w\rangle = \tau(\sum
 z_1 \langle z_2,w\rangle)=
 \tau(\sum S^2(w_1^\circ)\langle z,w_2\rangle)= \langle z, \sum
 \overline{\tau}(S^2(w_1^\circ))w_2\rangle= \langle z,
 ((\overline{\tau}S^2\circ)\star \id)(w)\rangle$.  The second, in the
 chain of
 equalities above, comes from Observation \ref{obse:equalities2}, part
 \,(4), equation \eqref{eqn:fourier2}.
 \end{proof}

We use the structure theorem --see Theorem \ref{theo:coalgcompacta}
and also Observation \ref{obse:equalities2}--
 to conclude that whenever $H$ is a compact quantum group, the maps
 $S_+: H \rightarrow H \,\, \text{and} \,\, S^2:H
 \rightarrow H$ are inner automorphisms.

\begin{theo}\label{theo:S2inner}
  Let $H$ be a compact quantum group
     with antipode $S$ and let $S_+$ be its positive antipode.
  Then, there exists a linear functional $\beta:H \rightarrow
    \mathbb{C}$ such that $S_+= (\beta \star \operatorname
    {id})(\operatorname {id} \star\, \beta^{-1})$ and  $S^2= (\beta^2 \star \operatorname
    {id})(\operatorname {id} \star\, \beta^{-2})$.

  Moreover, the restriction of $\beta$ to each simple
  component $H_\rho$, is unique up to multiplication by a non--zero
  scalar.

  Furthermore, we can assume that
  the linear functional $\beta$
   has been chosen as to satisfy   the following additional properties:
   \renewcommand{\theenumi}{\roman{enumi}}
   \begin{enumerate}
  \item  $\beta =\beta\, S_+$, $\beta^{-1}= \beta^{-1}\, S_+$;
  \item $\beta=\overline{\beta}\,\circ$;
   \item In the decomposition of the positive antipode $S_+=(\beta
     \star \id)(\id \star\, \beta^{-1})$, the operators $\beta \star \operatorname
    {id}$ and $\operatorname {id} \star\, \beta^{-1}$ are selfadjoint,
    commute with each other and leave invariat the simple components
     $H_\rho$ of $H$.
   \item The operators $\beta \star \operatorname {id}$ as well as
   $\operatorname {id} \star\, \beta^{-1}$ are positive.
   \end{enumerate}
   \end{theo}
   \begin{proof} We decompose $H=\bigoplus_\rho H_\rho$, where each
     $H_\rho$ is a simple coalgebra closed by the $\circ$
     operator of $H$. It was already proved that $S^2(H_\rho)=H_\rho$ and
   $S_+(H_\rho)=H_\rho$, hence
  for some $\beta_\rho:H_\rho \rightarrow \mathbb{C}$
   the automorphism $S_+|_{H_\rho}=(\beta_\rho \star \operatorname
   {id})(\operatorname {id} \star\, \beta^{-1}_\rho)$.
   If we call $\beta:H \rightarrow \mathbb{C}$
     the direct sum of all the $\beta_\rho$, we have that
  $S_+=(\beta \star \operatorname
   {id})(\operatorname {id} \star\, \beta^{-1})$.

  For the rest of the proof we proceed as follows.

   (i) Applying $\beta$ or $\beta^{-1}$ to the equality
$S_+(x)=\sum \beta(x_1) x_2
     \beta^{-1}(x_3)$ we deduce the result.

  (ii) Taking the adjoint of the operators in the
   formula $S_+=(\beta \star \id) (\id \star
     \beta^{-1})$ and recalling that $S_+$ is selfadjoint we
     obtain the following: $S_+=(\beta \star \id)^*(\id \star
     \beta^{-1})^*=
   ((\overline{\beta}S^2\circ) \star \id)(\id \star
   (\overline{\beta^{-1}}\circ))= ((\overline{\beta}\circ) \star\, \id)
   (\id \star (\overline{\beta}\circ)^{-1})$. Hence
  using again the uniqueness of $\beta$ up to a scalar in each simple
  coalgebra, we conclude that if we call $\beta_\rho$ the
  corresponding restriction, for each $\rho$ there exists a scalar
  $a_\rho \in \mathbb C$: $\beta_\rho= a_\rho
  \overline{\beta_\rho}\circ$. It is clear that $|a_\rho|=1$.
Indeed, we first deduce that $\beta_\rho \circ = a_\rho
\overline{\beta_\rho}$ and then --after conjugation and multiplication
 by $a_\rho$-- we obtain that $\beta_\rho=|a_\rho|^2 \beta_\rho$, i.e. $|a_\rho|=1$.

Hence if we take $b_\rho$ to be a solution of
  the equation $b_\rho a_\rho = \overline{b_\rho}$, and call
  $\gamma_\rho=b_\rho \beta_\rho$, the new functional $\gamma_\rho$
  satisfies the condition $\gamma_\rho = \overline{\gamma_\rho} \circ$.
  Hence, after changing $\beta_\rho$ by $\gamma_\rho$
whenever it is needed, we obtain a
  new functional $\beta$ --we call it $\beta$ instead of $\gamma$ to
  simplify notations-- that satisfies that $\beta =
  \overline{\beta}\, \circ$.

  (iii) The fact that the operators $\beta \star \id$ and $\id \star
   \beta^{-1}$ are selfadjoint is a direct consequence of (i) and
   (ii). Clearly they leave invariant the simple components $H_\rho$.

  (iv) This positivity result follows immediately from Corollary
  \ref{coro:albertpos}. Notice that the eventual change of sign
  necessary to make the operators positive --that we might have
  performed in accordance with the mentioned corollary-- does not affect
  the properties (i), (ii) and (iii).
  \end{proof}

Next we record for future use, some easy consequences of the results
about Radford's formula for coFrobenius coalgebras that appeared in
\cite{kn:bbt}. Previously, in \cite[Theorem 5.6]{kn:woro}
and in the context of
compact quantum groups, some of these topics were considered.

\begin{obse} \label{obse:nakayama} We say that $H$ is a coFrobenius
Hopf algebra over $\mathbb{C}$ if there exists a non--zero right (or
equivalently a non--zero left) integral that will be called $\varphi:H
\rightarrow  \mathbb{C}$. In this situation there are bijections
$x \mapsto (\varphi \leftharpoonup x) : H \rightarrow
H^{*,\,\text{rat}}$ and  $x \mapsto (x \rightharpoonup \varphi) : H \rightarrow
H^{*,\,\text{rat}}$.
\begin{enumerate}
\item Using the above isomorphisms, we can construct the Nakayama
automorphism of algebras characterized by the equality $\varphi(xy)=\varphi(y
 \mathcal N(x))$ --for all $x,y \in H$--. The modular
 function of the coFrobenius Hopf algebra
$H$, is defined as $\,\alpha= \varepsilon \mathcal N$.  The
distinghished group like element in $H$ --called $g$-- is
the unique group like element in $H$ such that for all $x \in H$: $\sum x_1 \varphi (x_2)= \varphi(x)g$. In
 \cite{kn:bbt} the authors prove the following formul\ae:
\begin{align} \label{eqn:nakayama} \mathcal N(x)= \sum \alpha(x_1)
  S^{-2}(x_2),\hspace*{.15cm}
S^4(x)=g\left(\sum \alpha(x_1)x_2 \alpha^{-1}(x_3)\right)g^{-1}
\hspace*{.15cm},\, x \in H,
\end{align}
where $g \in H$ is the distinguished group like element in $H$.

Clearly $\alpha:H \rightarrow \mathbb{C}$ is an algebra
  homomorphism and $\alpha^{-1}=
\alpha S$.
\item It follows easily from the equation \eqref{eqn:nakayama} that
  the inverse of $\mathcal N$ can be computed as:
\begin{equation}\label{eqn:nakayamainv}
\mathcal N^{-1}(x)= \sum \alpha^{-1}(x_1) S^2(x_2).
\end{equation}
Also, $\mathcal N$ satisfies the following formula:
 \begin{equation}\label{eqn:DeltaN}
\Delta(\mathcal N(x))= \sum \mathcal N(x_1) \otimes
S^{-2}(x_2)= \sum \mathcal N(x_1) \otimes
(\alpha^{-1} \star\, \mathcal N)(x_2).
\end{equation}
\noindent
In other words, the map
$\mathcal N$ is a morphism of comodules with the structures displayed
below:
$$\mathcal N :(H,\Delta)=H \rightarrow
(H,(\operatorname {id} \otimes\,  S^{-2})\Delta)= H_{S^{-2}}.$$

\item The Nakayama morphism is comultiplicative if and only if
  $\alpha=\varepsilon$ or equivalently if and only if $\mathcal N=
  S^{-2}$.
\noindent
Indeed, from the
comultiplicativity of $\mathcal N$ and from equation
\eqref{eqn:DeltaN} we deduce that for all $x \in H$, $\sum x_1
\otimes S^{-2}(x_2)=\sum x_1 \otimes \mathcal N(x_2)$. Then $\mathcal
N=S^{-2}$ and composing with $\varepsilon$ we deduce that
$\alpha=\varepsilon$.  The converse is clear.

\item If $H$ is also unimodular, we have that $\mathcal N^{-1}=
  S^{-1}\mathcal N S$.

\noindent
Indeed: \\
$\varphi(\mathcal N^{-1}(x) y) =
  \varphi(yx) = \varphi (Sx Sy)= \varphi (Sy \mathcal N Sx)= \varphi
  ((S^{-1}\mathcal N S )(x) y)$.
\item In the situation that the Hopf algebra is unimodular, the
  distinguished group like element $g=1$ and we have that:
\begin{equation}\label{eqn:radford}
 S^4 =  (\alpha \star \operatorname{id})(\operatorname{id} \star\,
 \alpha^{-1})
\end{equation}
\begin{equation}\label{eqn:radford2}
\mathcal N= \alpha \star S^{-2} = S^2 \star\, \alpha \hspace*{.25cm},\hspace*{.25cm}
\mathcal N^{-1}= \alpha^{-1}\star\, S^2 = S^{-2} \star\, \alpha^{-1}
\end{equation}
 It follows
  easily from the above formul\hspace*{-.02cm}\ae\, that
 $S^2 \mathcal N = \mathcal N S^2= \alpha
  \star \operatorname{id}$ and the
Nakayama automorphism commutes with $S^2$.
\end{enumerate}
\end{obse}

In the case of a compact quantum group more information can be
obtained about the maps considered above.

\begin{lema}\label{lema:nakayamapos} Let $H$ be a compact quantum
  group equipped with the inner product $\langle \,\,\,, \,\,\rangle$
given by the normal integral $\varphi$. Then:
\begin{enumerate}
\item The Nakayama
automorphism $\mathcal N$ is a positive operator with respect to
$\langle \,\,\,, \,\,\rangle$. Moreover it preserves the decomposition in simple
components $H_\rho$ --see Theorem \ref{theo:coalgcompacta}--.
\item The Nakayama automorphism commutes with $\circ$.
\item If $\alpha: H \rightarrow \mathbb C$ is the modular function,
  then $\alpha \star \operatorname {id} : H \rightarrow H$ is a
  positive operator with respect to the inner product $\langle \,\,\,, \,\,\rangle$.
\end{enumerate}
\end{lema}
\begin{proof} (1) The fact that $\mathcal N$ preserves the simple
  components $H_\rho$ follows immediately from Equation
  \eqref{eqn:nakayama}.
The following chain of equalities --applied for $y=x$--
proves the positivity of $\mathcal N$.
\begin{align} \label{eqn:Npositive} \langle \mathcal N x,y\rangle = \varphi(S(y^\circ)
  \mathcal N x)=\varphi(x S(y^\circ))=\langle S(y^\circ),S(x^\circ)\rangle.
\end{align}

For the last equality we used that: $S \circ S\, \circ =
\operatorname{id}$.
\noindent

(2) The proof of this commutation property follows from the equalities below:
 \begin{align} \label{eqn:Ncommutes} \varphi(y \mathcal
  N(x^\circ))= \varphi (x^\circ y)= \overline{\varphi (xy^\circ)} = \overline{\varphi
  (y^\circ \mathcal N(x))}=\varphi(y\mathcal N (x)^\circ).
\end{align}

\noindent

(3) This positivity result follows directly from the fact that
$\alpha \star\,\operatorname{id}$ is the composition of $\mathcal N$ and
$S^2$ that are positive and commuting operators.
\end{proof}

\begin{defi}\label{defi:posnak} Let $H$ be a compact quantum group and
call $\langle \,\,\,,\,\,\rangle$ the unitary inner product associated
to the normal integral $\varphi$. We define the positive square root
of Nakayama
automorphism --and denote it as $\mathcal P$--
as the unique automorphism of algebras that is
positive and satisfies that $\mathcal P^2=\mathcal N$.
\end{defi}

\begin{coro}\label{coro:squaremorf} Let $H$ be a compact quantum group
  and call $\varphi$ its normal integral. In the notations above, if
  $\beta$ is defined as in Theorem \ref{theo:S2inner}, then $\beta:H
  \rightarrow \mathbb C$ can be taken also to be an algebra morphism that
  satisfies
  $\beta^4=\alpha$.
\end{coro}
\begin{proof} The functional $\beta$ that we constructed before
satisfies: $S^4=\alpha \star \operatorname {id}
  \star \,\alpha^{-1}= \beta^4 \star \operatorname {id}\star\,
  \beta^{-4}$.
After restricting to each
  simple component $H_\rho$ we use Theorem \ref{theo:albert} to
prove the existence of a non zero scalar
  $c_\rho \in \mathbb{C}$
  , with the property that $\alpha|_{H_\rho} = c_\rho
  \beta^4|_{H_\rho}$. It is clear that being the operators associated
  to $\alpha|_{H_\rho}$ and to
$\beta|_{H_\rho}$ positive --see Lemma \ref{lema:nakayamapos} and
Theorem \ref{theo:S2inner}-- the scalar $c_\rho$ is in fact a positive
real number. If we change for each $\rho$ the functional
  $\beta_\rho$ by the functional $c_\rho^{1/4}\beta_\rho$ where
  $c_\rho^{1/4}$ is the {\em positive} fourth root of $ c_\rho$, we obtain
  that $\beta^4=\alpha$ and this change does not affect the positivity
  of $\beta$.

In order to prove the multiplicativity, we observe that in accordance
to the previous constructions $\beta \star \operatorname {id}$ is a
positive operator with the property that $(\beta \star \operatorname
{id})^4=\alpha \star \operatorname{id}$. It follows from Observation
\ref{obse:multiplicativity} that the unique fourth root of the
multiplicative map $\alpha\, \star \operatorname{id}$ is also
multiplicative. Hence, $\beta \star\,\operatorname {id}$ as well as
$\beta = \varepsilon (\beta \star\,\operatorname{id})$ are multiplicative.
\end{proof}

\begin{obse} \label{obse:importante}
\begin{enumerate}
\item Once we know that $\beta$ is
  a morphism of algebras it is clear that:
\begin{equation}\label{eqn:betaS}
\beta^{-1}= \beta S \hspace*{2cm} \beta=\beta^{-1}S.
\end{equation}
\item Below, we list the expressions of some of the maps considered
  above in terms of $\beta$.
\begin{equation} \label{eqn:Nbeta}
\mathcal N=\beta^2 \star \id \star \beta^2 \hspace*{.5cm} \mathcal
N^{-1}=\beta^{-2} \star \id \star \beta^{-2}
\end{equation}
\begin{equation*} 
\mathcal P=\beta \star \id \star \beta \hspace*{.5cm} S_+ =\beta \star \id \star \beta^{-1}
\end{equation*}
\begin{equation*} 
\mathcal P S_+ =\beta^2 \star \id  \hspace*{.5cm} \mathcal
P^{-1}S_+
=\id \star \beta^{-2}.
\end{equation*}

The formula for $\mathcal N$ follows directly from equation
\eqref{eqn:radford2}, and the formula for $\mathcal P$ follows from
the fact that $\id \star \beta$ and $\beta \star \id$ are positive
and from the uniqueness of the
square root of the positive operator $\mathcal N$.
\end{enumerate}
\end{obse}

Next we use the automorphism $\mathcal N$ to compute explicitly the
adjoint of $S$ and obtain some preliminary consequences of the normality of the
operator $S$. This initial result will be refined in Theorem
\ref{theo:meinanti}.
\begin{lema} \label{lema:S*} Let $H$ be a compact quantum group with
  normal
integral $\varphi$ and associated inner product $\langle\,\,\,,\,\,
  \rangle$.
\begin{enumerate}
\item The adjoint of the operator $S:H \rightarrow H$
  is $S^*=S\mathcal
  N^{-1} = \mathcal N S= S^{-1 }\star\, \alpha^{-1}$.
\item In the situation above, $SS^*= S^2 \mathcal N^{-1} = \mathcal N^{-1}S^2=
\operatorname{id} \star\, \alpha^{-1}$ and $S^*S= \mathcal N S^2= S^2
\mathcal N =
\alpha \star
\operatorname {id}$.
\item The operator $S$ is
  normal if and only if $S^4= \operatorname{id}$ and
  $\alpha^2=\varepsilon$.

\end{enumerate}
\end{lema}
\begin{proof} \begin{enumerate}
\item For $x,y \in H$ we have that:
$\langle Sx,y \rangle = \varphi((S\circ)y Sx)= \varphi(x y^\circ)$.
\noindent

Moreover:
$\langle x,S\mathcal N^{-1}y\rangle=\varphi
((S\!\circ\/\!S\mathcal N^{-1})(y)x)= \varphi((\circ \mathcal
N^{-1})(y)x)= \\
\overline{\varphi (\mathcal N^{-1}(y) x^\circ)}= \overline{\varphi
  (x^\circ y)}= \varphi(xy^\circ)$.

In the above chain of equalities we have used first that
$S\circ\/S=\circ$, then that $\varphi
\circ = \overline{\varphi}$, later
the definition of $\mathcal N$ to change the order of the product
inside $\varphi$ and then
again that $\varphi
\circ = \overline{\varphi}$, --see Observation \ref{obse:int1}--.
\item Clearly:
$SS^*=S^2 \mathcal N^{-1}= \mathcal N^{-1}S^2= (\alpha^{-1} \star
\operatorname {id})S^4=\operatorname {id} \star\, \alpha^{-1}$.

Also $S^*S= \mathcal N S^2= \alpha \star
\operatorname{id}$.

\item
Then, $S^*S= S^*S$ if and only if $\operatorname {id} \star\,\alpha^{-1} =
\alpha \star \operatorname{id}$.  Composing this equality  with
$\varepsilon$ we obtain that $\alpha^2=\varepsilon$. Hence, we deduce
that $\operatorname {id} \star\, \alpha= \alpha \star\,
\operatorname {id}$ and  this obviously implies that $S^4= (\alpha \star
\operatorname {id})(\operatorname {id} \star\, \alpha^{-1}) = \operatorname
{id}$. The converse follows easily.
\end{enumerate}
\end{proof}

\begin{obse} \label{obse:importante2}
Assume that we are in the situation considered above.
\begin{enumerate}
\item It is clear that the operators $SS^*$ and $S^*S$ leave all the
  simple coalgebras $H_\rho$ stable. Moreover the positive square
  roots of $SS^*$ and $S^*S$ are respectively: $\mathcal P^{-1}S_+$ and
  $\mathcal P S_+$. Then the right polar decomposition of $S$
  is $S=U \mathcal P S_+$ and the left polar
  decomposition is $S=\mathcal P^{-1}S_+U$, where
  $U$ is unitary and equal to: $U= SS_+^{-1}\mathcal
  P^{-1} = S_+^{-1}\mathcal PS$. It is clear that the
operator $U: H \rightarrow H$ is an antimultiplicative involution
--i.e. $U^2=\operatorname{id}$--. In the literature the
operator $U$ is sometimes called the {\em unitary antipode}.
\item In the case that $U =\operatorname{id}$, it is clear
that the antipode $S$, being equal to $S=\mathcal P S_+$ is
multiplicative. As it is also antimultiplicative and bijective, we
  conclude that the product of $H$ is commutative.
\end{enumerate}
\end{obse}

\begin{theo}\label{theo:meinanti} Let $H$ be a compact quantum
  group and consider the inner product given by the normal
  integral.
Then, the following properties are equivalent:
\begin{enumerate}
\item The compact quantum group $H$ is involutive: i.e., $S^2=\id$.
\item The positive antipode is trivial: i.e., $S_+= \operatorname {id}$.
\item The unitary antipode coincides with $S$: i.e., $U=S$.
\item The antipode $S$ has finite order.
\item The antipode $S$ is a normal operator: i.e., $S$ commutes with
  $S^*$.
\item The antipode $S$ is a selfadjoint operator.
\item The integral $\varphi:H \rightarrow \mathbb{C}$ is central:
  $\varphi(xy)=\varphi(yx)$ --see  \cite{kn:woro}--.
\item The automorphism $\mathcal N=\id$.
\end{enumerate}
\end{theo}
\begin{proof} As the only positive square root of the identity is the
  identity itself, (2) follows clearly from (1).
  If $S$ has finite order so does $S^2$ and also the positive
  automorphism $S_+$. Then, $S_+=\operatorname {id}$. This shows that
  the conditions (1), (2) and (4) are equivalent.
To prove that (1) implies (5) we observe that in the hypothesis of (1),
$\beta$ can be taken to be $\beta = \varepsilon$
and hence we deduce from Lemma \ref{lema:S*} that
$S$ is normal. Moreover, if we assume that the antipode is normal, using Lemma
\ref{lema:S*} we deduce that it has finite order.
We have proved the equivalence between (1), (2), (4) and (5). The
equivalence bewteen (6) and (7) is also clear. We know from Lemma
\ref{lema:S*} that $S^*= S\mathcal N^{-1}$. Hence $S$  is selfadjoint
if and only if $\mathcal N=\id$ and this happens --by definition of
$\mathcal N$-- if and only if $\varphi(xy)=\varphi(yx)$.

Now we prove that condition (2) implies
condition (6). Indeed in this case, it follows
that $\beta=\varepsilon$ and then --in accordance with equation
\eqref{eqn:Nbeta}-- we deduce that $\mathcal N=\id$.

We finish the proof by observing that condition (3) is equivalent to
condition (2). Indeed, in the case that $S=U$ as we know that $U^2 =
\operatorname {id}$, we conclude that $S^2=\operatorname
{id}$. Conversely, if $S_+=\operatorname {id}$, then condition (7)
guarantees that $\mathcal
N=\operatorname {id}$ and then $\mathcal P= S_+= \operatorname {id}$.
Using the formula $U=SS_+^{-1}\mathcal P^{-1}$ we conclude
that $U=S$ --see Observation
\ref{obse:importante2}--.

\end{proof}

Some additional informacion can be obtained about the unitary antipode
and Nakayama morphism in terms of the adjoint of the antipode.

\begin{theo}\label{theo:additional}
Let $H$ be a compact quantum
  group and consider the inner product given by the normal
  integral.
Then, the following properties are equivalent:
\begin{enumerate}
\item The adjoint of the antipode $S^*$ is an antimorphism of coalgebras.
\item The unitary antipode $U$ is a morphism of coalgebras.
\item The Nakayama morphism $\mathcal N$ is a comultiplicative morphism.
\item The functional $\alpha$ coincides with the counit.
\item The Nakayama morphism $\mathcal N$ equals $S^{-2}$.
\end{enumerate}
\end{theo}
\begin{proof}
Suppose (1) holds. Then $SS^*$ is a morphism of coalgebras
and so is $\mathcal PS_+$. Using that $S=U\mathcal PS_+$, we get that
$U$ is an antimorphism of coalgebras.

Suppose (2) holds. As $S=U\mathcal PS_+$, $S$ is an antimorphism of
coalgebras and $S_+$ is a morphism of coalgebras,
we deduce that $\mathcal P$ and $\mathcal N= \mathcal P^2$ are
morphisms of coalgebras and we have (3).

Suppose now that (3) holds. Then $\mathcal P$, $\mathcal PS_+$ and
$SS^*=(\mathcal PS_+)^2$ are comultiplicative.
Then, $S^*$ is an antimorphism of coalgebras.

Suppose now that (3) holds. From the equality $\mathcal N = \alpha
\star S^{-2}$ we conclude that $\alpha \star \operatorname {id}$ is
comultiplicative. In that case applying $\varepsilon \otimes
\varepsilon$ to the equality $\sum \alpha(x_1) x_2 \otimes \alpha(x_3)
x_4= \sum \alpha(x_1) x_2 \otimes x_3$ we deduce that $\alpha \star
\alpha = \alpha$ and then $\alpha = \varepsilon$, and we deduce that
condition (4) is satisfied.

It is clear that if $\alpha = \varepsilon$, then $\mathcal N = S^{-2}$
and then it is comultiplicative. In other words from condition (5) we
easily deduce condition (3).
\end{proof}
\section{The quantum special unitary group}
\label{section:su}

In this section we want to illustrate our constructions in the case of the
group $\su$, viewed as a $\circ$--Hopf algebra rather than as it is
usually presented in terms of a $\star$--operator.

\begin{defi} Let $\mu \in \mathbb{C}$ be a non zero real number such
  that $|\mu| < 1$. We call $\su=\mathbb{C}\langle \alpha, \gamma,
  \gamma^\circ, \widehat{\alpha}\rangle$ the non commutative algebra in
  the variables written above, subject to the following multiplicative
  relations.
\begin{equation}\label{eqn:msu}
\begin{split}
\widehat{\alpha}\alpha - \mu \gamma^\circ \gamma = 1&
\hspace*{1cm} \alpha \widehat{\alpha} - \mu^3 \gamma
\gamma^\circ = 1\\
\gamma^\circ & \gamma = \gamma \gamma^\circ\\
\mu \gamma \alpha = \alpha \gamma&\hspace*{1cm} \mu \gamma^\circ \alpha =
\alpha \gamma^\circ\\
\mu^{-1}\gamma \widehat{\alpha} = \widehat{\alpha}\gamma
&\hspace*{.5cm}\mu^{-1}\gamma^\circ \widehat{\alpha} =
\widehat{\alpha}\gamma^\circ
\end{split}
\end{equation} \label{sucirc}
The $\circ$ operator is defined as being multiplicative, involutive, conjugate
linear and taking on the generators $\alpha,\, \widehat{\alpha}$ the values:
\begin{equation}
\alpha^\circ = \alpha \hspace*{.5cm}
\widehat{\alpha}^\circ=\widehat{\alpha}
\end{equation}
The comultiplication and counit are defined as:

\begin{equation}\label{sucomult}
\begin{split} \Delta(\alpha)= \alpha \otimes \alpha + \mu^2
  \gamma^\circ \otimes \gamma \hspace*{1cm}& \Delta(\widehat{\alpha})= \widehat{\alpha} \otimes \widehat{\alpha} + \mu^2
  \gamma \otimes \gamma^\circ\\
\Delta({\gamma})=\gamma \otimes \alpha +\hspace*{1.2em} \widehat{\alpha}\hspace*{.3em} \otimes \gamma
\hspace*{1cm}& \Delta({\gamma^\circ})=\alpha \otimes \gamma^\circ +
\gamma^\circ \otimes \widehat{\alpha}
\end{split}
\end{equation}
\begin{equation}\label{eqn:suepsilon}
\varepsilon (\alpha)=1 \hspace*{1em} \varepsilon(\widehat{\alpha})=1
\hspace*{1em}  \varepsilon(\gamma)=0 \hspace*{1em} \varepsilon(\gamma^\circ)=0
\end{equation}
\end{defi}

The linear map $S: \su \rightarrow \su$ defined on the generators as:
\begin{equation}\label{eqn:suantipode}
S(\alpha)=\widehat{\alpha} \hspace*{1em} S(\widehat{\alpha})=\alpha
\hspace*{1em} S(\gamma)=-\mu\gamma \hspace*{1em} S(\gamma^\circ)=
-\mu^{-1}\gamma^\circ ,
\end{equation}
and then extended antimultiplicatively to all $\su$ is the antipode of
the above structure.

It is well known that $\su$ is a compact quantum group and as such it
is equipped with an integral and an inner product for which we shall
not provide an explicit expression. See \cite{kn:wororim} and
\cite{kn:woro} for a detailed  study of this family of Hopf algebras
from the $C^*$--algebra viewpoint, where explicit formul{\ae} are presented.

A direct computation shows that the linear functional $\theta:\su
\rightarrow \mathbb{C}$ defined by  the equality $S^2=\theta \star \id
\star \theta S$ can be taken as:
\begin{equation}\label{eqn:theta}
\theta(\alpha)=|\mu|^{-1} \hspace*{.5cm} \theta(\widehat{\alpha})=|\mu|
\hspace*{.5cm} \theta(\gamma)=\theta(\gamma^\circ)=0,
\end{equation}

and then extended multiplicatively.

In this situation, if we consider $\beta: \su \rightarrow \mathbb C$
defined as:
\begin{equation}\label{eqn:beta}
\beta(\alpha)=|\mu|^{-1/2} \hspace*{.5cm} \beta(\widehat{\alpha})=|\mu|^{1/2}
\hspace*{.5cm} \beta(\gamma)=\beta(\gamma^\circ)=0,
\end{equation}
we find that the positive antipode is given on the generators as:
\begin{equation}\label{eqn:s+su}
S_+(\alpha)=\alpha \hspace*{.5cm} S_+(\widehat{\alpha})=\widehat{\alpha}
\hspace*{.5cm} S_+(\gamma)=|\mu|\gamma \hspace*{.5cm}
S_+(\gamma^\circ)=|\mu|^{-1}\gamma^\circ.
\end{equation}
In this situation we obtain the following formul{\ae} for the Nakayama
automorphism $\mathcal N$ and for $\mathcal P$.

\begin{alignat}{4}\label{eqn:nakayamabeta}
\mathcal N(\alpha)&=\mu^{-2}\alpha &\hspace*{.5cm} \mathcal N(\widehat{\alpha})&= \mu^2 \widehat{\alpha}
&\hspace*{.5cm} \mathcal
N(\gamma)&= \gamma &\hspace*{.5cm} \mathcal N(\gamma^\circ)&=
\gamma^\circ\\
\mathcal P(\alpha)&=|\mu|^{-1}\alpha &\hspace*{.5cm} \mathcal P(\widehat{\alpha})&= |\mu|\widehat{\alpha}
&\hspace*{.5cm} \mathcal
P(\gamma)&= \gamma &\hspace*{.5cm} \mathcal P(\gamma^\circ)&=
\gamma^\circ
\end{alignat}

The unitary antipode $U$ can be computed by the following formulae:
\begin{align}\label{eqn:unitaryantipodebeta}
U(\alpha)=|\mu|\widehat{\alpha} \hspace*{.5cm}U(\widehat{\alpha})=
|\mu|^{-1} \alpha \hspace*{.5cm}
U(\gamma)=-\operatorname{sg}(\mu)\gamma \hspace*{.5cm} U(\gamma^\circ)=-\operatorname{sg}(\mu)\gamma^\circ
\end{align}

\end{document}